\documentclass[12pt]{article}
\usepackage{amsmath}
\usepackage{amsfonts}
\usepackage{amsthm}
\usepackage{graphics}
\numberwithin{equation}{section}
\renewcommand{\O}{\mathcal{O}}
\newtheorem{definition}{Definition}[section]

\newtheorem{lemma}{Lemma}[section]
\begin{document}
\title{Discrete Painlev\'e equations for recurrence coefficients of orthogonal polynomials}         
\author{Walter Van Assche \thanks{Research supported by INTAS Research Network NeCCA (03-51-6637), FWO projects 
G.0184.02 and G.0455.04, and OT/04/21 of the Research Counsel of K.U.Leuven.}\\
Katholieke Universiteit Leuven, Belgium}        
\date{}          
\maketitle

\section{Introduction}       
Orthonormal polynomials on the real line are defined by the orthogonality conditions
\begin{equation}  \label{eq:orthonormal}
  \int_{\mathbb R} p_n(x)p_m(x)\, d\mu(x) = \delta_{m,n}, \qquad m,n \geq 0,
\end{equation}
where $\mu$ is a positive measure on the real line for which all the moments exist and
$p_n(x) = \gamma_n x^n + \cdots$, with positive leading coefficient $\gamma_n > 0$.
A family of orthonormal polynomials always satisfies a three-term recurrence relation
of the form
\begin{equation}  \label{eq:3tr}
   xp_n(x) = a_{n+1}p_{n+1}(x) + b_n p_n(x) + a_n p_{n-1}(x), \qquad n \geq 0,
\end{equation}
with $p_{-1}=0$ and 
\begin{equation} \label{eq:gamma0}  
p_0=\gamma_0 = \left(\int d\mu(x) \right)^{-1/2}. 
\end{equation}
Comparing the leading coefficients in (\ref{eq:3tr}) gives
\begin{equation}  \label{eq:angamma}
   a_{n+1} = \frac{\gamma_n}{\gamma_{n+1}} > 0,
\end{equation}
and computing the Fourier coefficients of $xp_n(x)$ in (\ref{eq:3tr}) gives
\begin{eqnarray}
    a_n & = & \int_{\mathbb R} xp_n(x)p_{n-1}(x)\, d\mu(x),  \label{eq:an} \\
    b_n & = & \int_{\mathbb R} xp_n^2(x)\, d\mu(x).  \label{eq:bn}
\end{eqnarray}
The converse statement is also true and is known as \textit{the spectral theorem for orthogonal polynomials}:
if a family of polynomials satisfies a three-term recurrence relation of the form (\ref{eq:3tr}),
with $a_n > 0$ and $b_n \in \mathbb{R}$ and initial conditions $p_0=1$ and $p_{-1}=0$, then there
exists a \textit{probability measure} $\mu$ on the real line such that these polynomials are
orthonormal polynomials satisfying (\ref{eq:orthonormal}).
This gives rise to two important problems:
\begin{description}
\item[Problem 1.] Suppose the measure $\mu$ is known. What can be said about the recurrence coefficients
$\{a_n: n =1,2,3,\ldots\}$ and $\{b_n: n=0,1,2,\ldots\}$? This is known as the \textit{direct problem} for
orthogonal polynomials.
\item[Problem 2.] Suppose the recurrence coefficients $\{a_{n+1},b_n: n=0,1,2,\ldots\}$ are known. What can
be said about the orthogonality measure $\mu$? This is known as the \textit{inverse problem} for orthogonal polynomials.
\end{description}
The recurrence coefficients are usually collected in a tridiagonal matrix of the form
\begin{equation}  \label{eq:jacobi}
    J = \begin{pmatrix}
           b_0 & a_1 & & & \\
           a_1 & b_1 & a_2 & & \\
               & a_2 & b_2 & a_3 & \\
               &     & a_3 & b_3 & \ddots \\
	       &     &     & \ddots & \ddots
          \end{pmatrix}  
\end{equation} 	
which acts as an operator on (a subset of) $\ell_2({\mathbb N})$ and which is known as the Jacobi matrix or
Jacobi operator. If $J$ is self-adjoint, then the spectral measure for $J$ is precisely the orthogonality
measure $\mu$. Hence problem 1 corresponds to the inverse problem for the Jacobi matrix $J$ and problem 2
corresponds to the direct problem for $J$.

In the present paper we will study problem 1 for a few special cases. In Section 2 we study measures on the
real line with an exponential weight function of the form $d\mu(x) = |x|^\rho \exp(-|x|^m)\, dx$, which are
known as Freud weights, named after G\'eza Freud who studied them in the 1970's. It will be shown that the
recurrence coefficients $a_n$ satisfy a non-linear recurrence relation which corresponds to the discrete
Painlev\'e I equation and its hierarchy. In Section 3 we will study a family of orthogonal polynomials on the
unit circle. We will first give some background on orthogonal polynomials on the unit circle and the
corresponding recurrence relations. We will study the weight function $w(\theta) = \exp(\lambda \cos \theta)$, 
and it will be shown that the recurrence coefficients satisfy a non-linear recurrence relation which corresponds
to discrete Painlev\'e II. These orthogonal polynomials play an important role in random unitary matrices
and combinatorial problems for random permutations. In Section 4 we will study certain discrete
orthogonal polynomials related to Charlier polynomials. The recurrence coefficients $a_n$ and $b_n$ are shown
to satisfy a system of non-linear recurrence relations which are again related to the discrete Painlev\'e II equation.
Finally, in Section 5 we consider certain $q$-orthogonal polynomials which are $q$-analogs of the Freud weight.
We will show that the recurrence coefficients satisfy a $q$-deformed Painlev\'e I equation.
Most of the material in this paper is not new: the recurrence relations in Section 2 were already obtained
by Freud in \cite{freud} and are known as \textit{Freud equations} in the field of orthogonal polynomials.
It was Magnus \cite{magnus2} who made the connection with the discrete Painlev\'e equation I.
The recurrence relation in Section 3 was found by Periwal and Shevitz \cite{periwal} (see also
Hisakado \cite{hisakado},  Tracy and Widom \cite{tracy}; Baik \cite{baik0} used the Riemann-Hilbert
approach to obtain the Painlev\'e equation). The recurrence relations in Section 4 were obtained
by Van Assche and Foupouagnigni in \cite{wvafoup}. The results in Section 5 were first obtained by Nijhoff
\cite{nijhoff}. We hope that bringing together these various orthogonal
polynomials and the corresponding discrete Painlev\'e equations will be illuminating and encourage researchers
in the field of orthogonal polynomials and researchers in integrable difference equations 
to talk to each other, and that the interaction between both areas of mathematics will shed some extra light on
either subject. 
 
\section{Freud weights}   
Freud weights are exponential weights on the real line $(-\infty,\infty)$ of the form
\[   w_\rho(x) = |x|^\rho \exp(-|x|^m), \qquad \rho > -1, m > 0. \]
They were considered by G\'eza Freud in his 1976 paper \cite{freud}, where he gave
a recurrence relation for the recurrence coefficients $a_n$ when $m=2,4,6$. For these cases
Freud found the asymptotic behavior of the recurrence coefficients $a_n$ and he formulated
a conjecture for this asymptotic behavior for every $m > 0$. This conjecture really started
the analysis of general orthogonal polynomials on unbounded sets of the real line, since before
Freud's work only very special cases such as the Hermite polynomials and the Laguerre polynomials
were studied in detail. In this section I will repeat Freud's analysis of the recurrence coefficients
of Freud weights, make the connection with discrete Painlev\'e equations (which was not known to
Freud but first pointed out by Magnus in \cite{magnus2}), and point out what has been done after Freud. 
Observe that Freud weights are symmetric, i.e.,  $w_\rho(-x)=w_\rho(x)$, which implies that $b_n=0$ for $n\geq 0$.

\subsection{Generalized Hermite polynomials}
The case $m=2$ corresponds to generalized Hermite polynomials (and $\rho=0$ are the Hermite polynomials).
Generalized Hermite polynomials were already investigated by Chihara in \cite{chihara1} (see also \cite{chihara}).
The weight $w_0(x) = \exp(-x^2)$ satisfies the first order differential equation
\begin{equation}  \label{eq:pearson2}
   [w_0(x)]' = -2x w_0(x) 
\end{equation}
which is the \textit{Pearson equation} for the Hermite weight. In general a weight $w$ satisfying a Pearson equation
of the form $[\sigma(x)w(x)]' = \tau(x) w(x)$, with $\sigma$ a polynomial of degree at most two and $\tau$ a polynomial
of degree one, is called a classical weight. The weight functions for Hermite polynomials, Laguerre polynomials,
and Jacobi polynomials are the classical weights for $\sigma$ of degree zero, one and two, respectively. Bessel polynomials
appear for $\sigma(x)=x^2$, but they are not orthogonal on the real line with respect to a positive measure.
Freud's idea was to compute the integral
\begin{equation}  \label{eq:mainint}
   \int_{-\infty}^\infty \left( p_n(x)p_{n-1}(x) |x|^\rho \right)' w_0(x)\, dx 
\end{equation}
in two different ways. The first way is simply working out the derivative in the integrand and to use
the orthogonality to evaluate the resulting terms. This gives
\begin{eqnarray}
  \int_{-\infty}^\infty \left( p_n(x)p_{n-1}(x) |x|^\rho \right)' w_0(x)\, dx & = &
   \int_{-\infty}^\infty  p_n'(x)p_{n-1}(x) |x|^\rho  w_0(x)\, dx   \label{eq:int1} \\
  & & +\ \int_{-\infty}^\infty p_n(x)p_{n-1}'(x) |x|^\rho  w_0(x)\, dx \label{eq:int2} \\
  & & +\ \rho \int_{-\infty}^\infty \frac{p_n(x)p_{n-1}(x)}{x} |x|^\rho  w_0(x)\, dx.  \label{eq:int3} 
\end{eqnarray}  
For the right hand side in (\ref{eq:int1}) we use the fact that
\begin{eqnarray*}    
  p_n'(x) &=& n \gamma_n x^{n-1} + \textrm{ lower order terms } \\
          &=& n \frac{\gamma_n}{\gamma_{n-1}} p_{n-1}(x)      + \textrm{ lower degree terms}. 
\end{eqnarray*}
This gives
\[  \int_{-\infty}^\infty  p_n'(x)p_{n-1}(x) |x|^\rho  w_0(x)\, dx = n \frac{\gamma_n}{\gamma_{n-1}} . \]
For the integral in (\ref{eq:int2}) we see that $p_{n-1}'(x)$ is a polynomial of degree $n-2$ and hence by orthogonality the integral vanishes. For the integral in (\ref{eq:int3}) we use the fact that the weight function $w_\rho$ is even, i.e.,
$w_\rho(-x)=w_\rho(x)$, which implies that $p_n(-x)=(-1)^n p_n(x)$. This means that
$p_n(x)/x$ is a polynomial of degree $n-1$ when $n$ is odd and $p_{n-1}(x)/x$ is a polynomial of degree $n-2$ when
$n$ is even. Hence when $n$ is even the integral in (\ref{eq:int3}) vanishes, and when $n$ is odd we have
\[   \frac{p_n(x)}{x} = \frac{\gamma_n}{\gamma_{n-1}} p_{n-1}(x) + \textrm{ lower degree terms}, \]
so that
\[   \rho \int_{-\infty}^\infty \frac{p_n(x)p_{n-1}(x)}{x} |x|^\rho  w_0(x)\, dx = \rho \frac{\gamma_n}{\gamma_{n-1}} \Delta_n \]
where 
\[  \Delta_n = \begin{cases} 0 & \textrm{if $n$ is even,} \\
                              1 &   \textrm{if $n$ is odd.}
               \end{cases}  \]
Combining these results and using the expression $a_n = \gamma_{n-1}/\gamma_n$ gives
\begin{equation}  \label{eq:mainint1}
   \int_{-\infty}^\infty \left( p_n(x)p_{n-1}(x) |x|^\rho \right)' w_0(x)\, dx 
   = \frac{n+\rho \Delta_n}{a_n}. 
\end{equation}
Observe that this holds whenever $w$ is a symmetric weight on the real line.
A second way to compute the integral in (\ref{eq:mainint}) is to use integration by parts, combined with
Pearson's equation (\ref{eq:pearson2}) for the weight. This gives
\begin{eqnarray}  \int_{-\infty}^\infty \left( p_n(x)p_{n-1}(x) |x|^\rho \right)' w_0(x)\, dx
  & = & - \int_{-\infty}^\infty  p_n(x)p_{n-1}(x) |x|^\rho w_0'(x)\, dx  \nonumber \\
  & = & 2 \int_{-\infty}^\infty  xp_n(x)p_{n-1}(x) |x|^\rho  w_0(x)\, dx \nonumber \\
  & = & 2 a_n .  \label{eq:twoan}
\end{eqnarray}
Combining (\ref{eq:mainint1}) and (\ref{eq:twoan}) then gives
\begin{equation}  \label{eq:anH}
              a_n^2 = \frac{n+\rho \Delta_n}{2}
\end{equation}
so that $a_n = \sqrt{n+\rho \Delta_n}/\sqrt{2}$. For $\rho=0$, which are the Hermite polynomials, one has $a_n=\sqrt{n/2}$.
For generalized Hermite polynomials one has
\begin{equation}  \label{eq:limit2}
  \lim_{n \to \infty} \frac{a_n}{\sqrt{n}} = \frac{1}{\sqrt{2}}.
\end{equation}

\subsection{Freud weight \boldmath $m=4$}
Let us now consider the weight $w_\rho(x) = |x|^\rho \exp(-x^4)$, where $\rho > -1$. 
The Pearson equation for the weight $w_0(x)=\exp(-x^4)$ is
\[   [w_0(x)]' = -4 x^3 w_0(x), \]
which is a first order linear differential equation with polynomial coefficients, but since the polynomial
coefficient is now of degree $3$, the weight $w_0$ is no longer a classical weight but a semi-classical one.
The equation (\ref{eq:mainint1}) remains valid for this Freud weight, but integration by parts give a different
result. Indeed
\begin{eqnarray}  \int_{-\infty}^\infty \left( p_n(x)p_{n-1}(x) |x|^\rho \right)' w_0(x)\, dx
  & = & - \int_{-\infty}^\infty  p_n(x)p_{n-1}(x) |x|^\rho w_0'(x)\, dx  \nonumber \\
  & = & 4 \int_{-\infty}^\infty  x^3p_n(x)p_{n-1}(x) |x|^\rho  w_0(x)\, dx.  \label{eq:part}
\end{eqnarray}
This integral can be computed by applying the three term recurrence (\ref{eq:3tr}), with $b_n=0$, repeatedly. Indeed
\begin{eqnarray*}
   x^2 p_n(x) &=& a_{n+1} x p_{n+1}(x) + a_n x p_{n-1}(x) \\
         & = &  a_{n+1}a_{n+2} p_{n+2}(x) + (a_{n+1}^2+a_n^2) p_{n}(x) + a_na_{n-1} p_{n-2}(x),
\end{eqnarray*}
and   
\begin{eqnarray*}
   x^3 p_n(x)  & = &  a_{n+1}a_{n+2} xp_{n+2}(x) + (a_{n+1}^2+a_n^2) xp_{n}(x) + a_na_{n-1} xp_{n-2}(x) \\
     & = & a_{n+1}a_{n+2}a_{n+3} p_{n+3}(x) + a_{n+1}(a_{n+2}^2+a_{n+1}^2+a_n^2) p_{n+1}(x) \\
    & & +\ a_n(a_{n+1}^2+a_n^2+a_{n-1}^2) p_{n-1}(x) + a_na_{n-1}a_{n-2} p_{n-3}(x).
\end{eqnarray*}
 From this one easily finds
\begin{equation} \label{eq:x3}
  \int_{-\infty}^\infty x^3 p_n(x)p_{n-1}(x) |x|^\rho w_0(x)\, dx = a_n(a_{n+1}^2+a_n^2+a_{n-1}^2).
\end{equation}
This result holds for $n \geq 1$ if we define $a_0=0$.
Observe that one can obtain this also by using the calculus of the Jacobi operator, since
\[   (J^3)_{m,n} = \int_{-\infty}^\infty x^3 p_m(x)p_n(x)w_\rho(x)\, dx, \]
and the quantity of interest is $(J^3)_{n,n-1}$. This computation is quite simple since it amounts to
some simple matrix multiplications. Combining (\ref{eq:mainint1}) and  (\ref{eq:part})--(\ref{eq:x3}) then gives
\begin{equation} \label{eq:freud4}
     4a_n^2(a_{n+1}^2 + a_n^2 + a_{n-1}^2) = n + \rho \Delta_n .
\end{equation}
This time we do not get $a_n$ explicitly, but instead we get a second order non-linear recurrence
relation for the recurrence coefficients. The initial conditions are $a_0=0$ and for $a_1$ we require that
$p_1(x) = xp_0(x)/a_1$ has norm one, which means
\[     \frac{\gamma_0^2}{a_1^2} \int_{-\infty}^\infty x^2 w_\rho(x)\, dx = 1, \]
and together with (\ref{eq:gamma0}) this means
\[   a_1 = \left( \frac{\Gamma(\frac{3+\rho}4)}{\Gamma(\frac{1+\rho}{4})} \right)^{1/2} . \]
If we put $x_n = 2 a_n^2$ then clearly $x_0=0$, $x_n > 0$ for $n > 0$ and
\begin{equation}  \label{eq:Freud4x}
     x_n(x_{n+1}+x_n+x_{n-1}) = n + \rho \Delta_n .
\end{equation}
This recurrence relation is the \textit{discrete Painlev\'e equation} d-P$_{\textrm I}$
\[   x_{n+1} + x_n + x_{n-1} = \frac{z_n + \gamma (-1)^n}{x_n} + \delta, \qquad z_n = \alpha n + \beta, \]
with $\alpha=1$, $\beta=\rho/2$, $\gamma=-\rho/2$ and $\delta=0$, since we can write $\Delta_n=(1-(-1)^n)/2$.

The observation that Freud's equation (\ref{eq:freud4}) is a discrete Painlev\'e equation was not
known to Freud but was pointed out much later by Magnus in \cite{magnus2}. This means that the equation
has the discrete Painlev\'e property, which is known as \textit{singularity confinement} (see, e.g., \cite{GNR}):
\begin{definition}[discrete Painlev\'e property]
If $x_n$ is such that it results in a singularity for $x_{n+1}$, then there exists a $p \in \mathbb{N}$ such that
this singularity is confined to $x_{n+1},\ldots,x_{n+p}$. Furthermore $x_{n+p+1}$ depends only
on $x_{n-1},x_{n-2},\ldots$.
\end{definition}
The usual Painlev\'e property for differential equations is that the only movable singularities (singularities
which depend on the initial conditions) of solutions
of a Painlev\'e equation are poles. Poles are isolated singularities, hence a discrete version of poles
as singularities is to require that singularities of a discrete equation are confined. This is the case for discrete Painlev\'e I. Consider for instance d-P$_{\textrm I}$ in the form
\[     x_n(x_{n+1}+x_n+x_{n-1}) = n, \]
then 
\[     x_{n+1} = \frac{n}{x_n} - x_{n} - x_{n-1}, \]
and if $x_n=0$ then we have a singularity for $x_{n+1}$ which becomes $\pm \infty$. For $x_{n+2}$ we then find
$\mp \infty$ and for $x_{n+3}$ we have the indeterminate form $(\pm \infty) + (\mp \infty)$. A more careful
analysis is to put $x_n = \epsilon$ and to expand $x_{n+k}$ in a Laurent series in $\epsilon$. This gives
\begin{eqnarray*}
   x_n & = & \epsilon \\[8pt]
   x_{n+1} & = & \frac{n}{\epsilon} - x_{n-1} - \epsilon \\
   x_{n+2} & = & - \frac{n}{\epsilon} + x_{n-1} + \frac{n+1}{n} \epsilon + \O(\epsilon^2) \\ 
   x_{n+3} & = & - \frac{n+3}{n} \epsilon + \O(\epsilon^2) \\
   x_{n+4} & = & \frac{n}{n+3} x_{n-1} + \O(\epsilon).
\end{eqnarray*}
So we see that as $\epsilon \to 0$ the indeterminate form for $x_{n+3}$ becomes 0, but it does not give
a new singularity for $x_{n+4}$. The singularity is confined to $x_{n+1},x_{n+2},x_{n+3}$.

The solution of d-P$_{\textrm I}$ can not be obtained in a closed form, but one can say a few things about the behavior
of the solution. Freud obtained the asymptotic behavior of the solution of (\ref{eq:Freud4x}) in the following way. 
Since $x_n \geq 0$ for $n \geq 0$ we have
\[   x_n^2 \leq x_n(x_{n+1}+x_n+x_{n-1}) = n + \rho \Delta_n, \]
so that $x_n/\sqrt{n}$ is a bounded and positive sequence. Define $A$ and $B$ as the smallest and largest
accumulation points
\[    0 \leq A = \liminf_{n \to \infty} \frac{x_n}{\sqrt{n}} 
\leq \limsup_{n \to \infty} \frac{x_n}{\sqrt{n}} = B < \infty . \]
Choose a subsequence $n'$ such that $x_{n'}/\sqrt{n'} \to A$, then as $n' \to \infty$ in (\ref{eq:Freud4x}) we
have
\[   1 \leq A(2B+A). \]
In a similar way we can choose a subsequence $n''$ such that $x_{n''}/\sqrt{n''} \to B$, and as $n' \to \infty$ in (\ref{eq:Freud4x}) we then find
\[    B(2A+B) \leq 1. \]
Together this gives
\[     B(2A+B) \leq A(2B + A), \]
so that $B^2 \leq A^2$. But since we already know that $A \leq B$, this implies that $A=B$ and hence
\[  \lim_{n \to \infty} \frac{x_n}{\sqrt{n}} = A = B . \]
If we take the limit in the recurrence relation (\ref{eq:Freud4x}) then one finds $3A^2=1$ so that $A=1/\sqrt{3}$.
Recall that $x_n=2a_n^2$, hence as a result we have
\begin{equation}  \label{eq:limit4}
  \lim_{n \to \infty} \frac{a_n}{n^{1/4}} = \frac{1}{\sqrt[4]{12}} . 
\end{equation}

\begin{figure}[ht]
\centering
\rotatebox{-90}{\resizebox{3.5in}{!}{\includegraphics{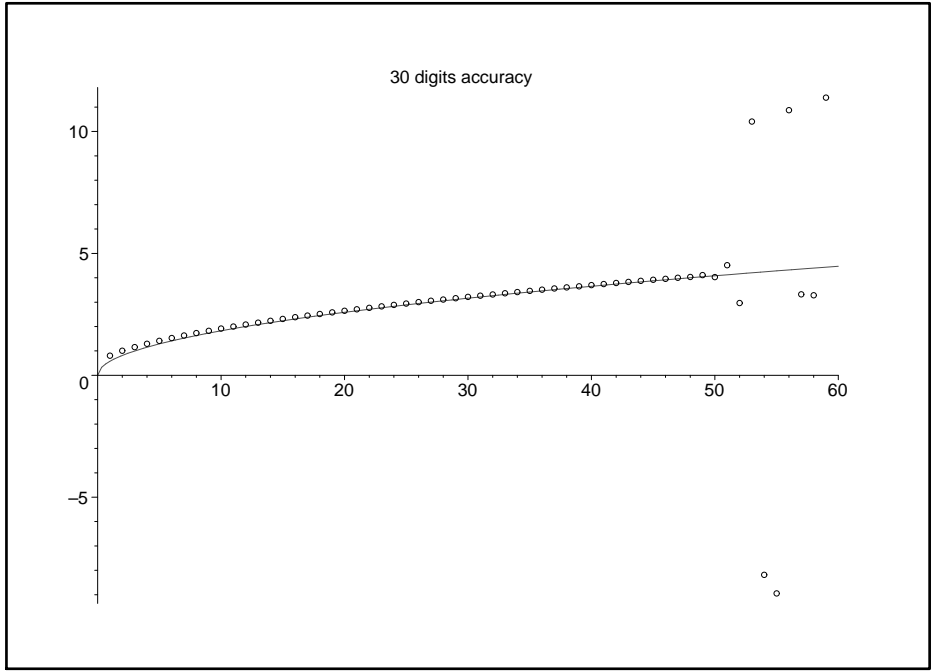}}}
\caption{The result of computing $x_n$ from d-P$_{\textrm I}$ using 30 significant digits}
\label{fig:freud4}
\end{figure}

The recurrence relation (\ref{eq:Freud4x}), with initial conditions
\[   x_0 = 0, \quad x_1 = \frac{2\Gamma(\frac{3+\rho}{4})}{\Gamma(\frac{1+\rho}{4})} \]
 is very unstable for computing the recurrence coefficients.  Lew and Quarles \cite{lew} showed that there is
a unique positive solution of the recurrence relation (\ref{eq:Freud4x}) with $x_0=0$. Hence a small
error in $x_1$ eventually destroys the positivity of $x_n$. In Figure \ref{fig:freud4} we plotted the values
of $x_n$ obtained from the recurrence relation (with $\rho = 0$) by using 30 digits accuracy. The $x_n$
are following the asymptotic behavior $\sqrt{n/3}$ quite well until $n \sim 50$ and then large deviations from the true solution appear. 

Lew and Quarles proved the unicity by showing that there is an operator $F$ acting on a Banach space of
infinite sequences with $x_0=0$, such that the positive solution $x=(x_0,x_1,x_2,\ldots)$ of (\ref{eq:Freud4x}) 
is a fixed point: $F(x)=x$. The operator $F$ is then shown to be a contraction, so that the fixed point is unique.

Observe that if we take a weight function of the form $w_\rho(x) = |x|^\rho \exp(-x^4 + \lambda x^2)$, then the
Pearson equation becomes
\[   [w_0(x)]' = (-4x^3+2\lambda x) w_0(x), \]
and a slight modification of the previous computations gives the recurrence
\[   4a_n^2 (a_{n+1}^2 + a_n^2 + a_{n-1}^2) - 2\lambda a_n^2 = n+ \rho \Delta_n. \]
If we put $x_n=2a_n^2$, then the $x_n$ satisfy the discrete Painlev\'e equation d-P$_{\textrm I}$ with
$\alpha=1$, $\beta=\rho/2$, $\gamma=-\rho/2$ and $\delta=\lambda$.

\subsection{Freud weight \boldmath $m=6$}
For the Freud weight $w_\rho(x) = |x|^\rho \exp(-x^6)$ one can proceed in a very similar way. The Pearson
equation now becomes
\[   [w_0(x)]' = -6 x^5 w_0(x), \]
so that
\[ \int_{-\infty}^\infty \left( p_n(x)p_{n-1}(x) |x|^\rho \right)' w_0(x)\, dx
  = 6 \int_{-\infty}^\infty  x^5 p_n(x)p_{n-1}(x) |x|^\rho  w_0(x)\, dx. \]
The integral on the right is 
\[   \int_{-\infty}^\infty  x^5 p_n(x)p_{n-1}(x) |x|^\rho  w_0(x)\, dx = (J^5)_{n,n-1} \]
and this is
\begin{multline*}
  (J^5)_{n,n-1} = a_n (a_{n-2}^2a_{n-1}^2 + a_{n-1}^4 + 2 a_{n-1}^2a_{n}^2 + a_{n-1}^2a_{n+1}^2 \\
 +\ a_n^4 + 2 a_n^2a_{n+1}^2 + a_{n+1}^4 + a_{n+1}^2 a_{n+2}^2). 
\end{multline*}
Together with (\ref{eq:mainint1}) this gives
\begin{multline}  \label{eq:freud6}
  6a_n^2 (a_{n-2}^2a_{n-1}^2 + a_{n-1}^4 + 2 a_{n-1}^2a_{n}^2 + a_{n-1}^2a_{n+1}^2 \\
 +\ a_n^4 + 2 a_n^2a_{n+1}^2 + a_{n+1}^4 + a_{n+1}^2 a_{n+2}^2) = n+\rho \Delta_n. 
\end{multline} 
This is a fourth order non-linear recurrence relation for the recurrence coefficients. It is within the hierarchy of
discrete Painlev\'e I \cite{joshi}. Freud was also able to obtain the asymptotic behavior for the $a_n$ in this case.
Obviously
\[  6a_n^6 \leq  6a_n^2 (a_{n-2}^2a_{n-1}^2 + a_{n-1}^4 + 2 a_{n-1}^2a_{n}^2 + a_{n-1}^2a_{n+1}^2
 + a_n^4 + 2 a_n^2a_{n+1}^2 + a_{n+1}^4 + a_{n+1}^2 a_{n+2}^2) = n+\rho \Delta_n, \]
hence $a_n/n^{1/6}$ is a positive and bounded sequence. If we define
\[   0 \leq A = \liminf_{n \to \infty} \frac{a_n^2}{n^{1/3}} 
\leq \limsup_{n \to \infty} \frac{a_n^2}{n^{1/3}} = B < \infty ,\]
then by taking a subsequence that converges to $A$ we find
\[ 6A(5B^2 + 4AB + A) \geq 1 \]
and by taking a subsequence that converges to $B$ we find
\[  6B(5A^2 + 4AB + B^2) \leq 1. \]
Combining both inequalities gives
\[  6B(5A^2 + 4AB + B^2) \leq 6A(5B^2 + 4AB + A) \]
so that $A^2B + B^3 \leq AB^2+A^3$. This is equivalent with
\[  (B-A)(A^2+AB+B^2) \leq (B-A) AB. \]
If we assume that $A<B$, then this would imply that $A^2+B^2 \leq 0$, or $A=B=0$, which is impossible (since $A<B$). Hence
our assumption is false and $A=B$. Taking the limit in (\ref{eq:freud6}) gives $60A^3 = 1$, hence
\begin{equation}  \label{eq:limit6}
  \lim_{n \to \infty} \frac{a_n}{n^{1/6}} = \frac{1}{\sqrt[6]{60}}.
\end{equation}

\subsection{Freud's conjecture} 
On the basis of (\ref{eq:limit2}), (\ref{eq:limit4}) and (\ref{eq:limit6}), Freud made the conjecture that
for every $m > 0$ and $\rho > -1$ one has
\begin{equation} \label{eq:freudcon}
  \lim_{n \to \infty} \frac{a_n}{n^{1/m}} 
= \left( \frac{\Gamma(m+1)}{\Gamma(\frac{m}2)\Gamma(\frac{m}2 + 1)} \right)^{-1/m}.  
\end{equation}
Furthermore Freud showed that if the limit exists for even $m$, then it is equal to the expression in (\ref{eq:freudcon}).
Freud could not prove the existence of the limit for even $m \geq 8$ because the central coefficient $a_n$
in the non-linear recurrence relation does not occur sufficiently often and more non central coefficients
$a_{n\pm k}$ with $k\neq 0$ appear, making the recurrence no longer `diagonally dominant'. The simple trick
using $\limsup$ and $\liminf$ then no longer suffices to show that the limit exists. The proof of Freud's conjecture
for every even integer $m$ was given by Alphonse Magnus \cite{magnus} \cite{magnus1}. His proof still consists
of obtaining a non-linear recurrence relation for the $a_n$ (the Freud equation, which is within the hierarchy of 
discrete Painlev\'e I), but a more subtle argument is used to prove the existence of the limit.
Freud's conjecture for general $m >0$ was finally proved by Lubinsky, Mhaskar and Saff \cite{lubinsky}. The proof
for general $m > 0$ no longer uses a recurrence relation for the recurrence coefficients but relies
on the Mhaskar-Rakhmanov-Saff number and results of weighted polynomial approximation and an equilibrium
problem of logarithmic potential theory with external field.  

For $m$ an even positive integer, M\'at\'e, Nevai and Zaslavsky \cite{mate} obtained an asymptotic expansion of the form
\[   \frac{a_n^2}{n^{1/m}} = \sum_{k=0}^{\infty}   \frac{c_k}{n^{2k}}, \]
where $c_0$ is the constant in Freud's conjecture (\ref{eq:freudcon}), 
but the other coefficients $c_k$ with $k > 0$ are not explicitly known. Their analysis is again based on the non-linear recurrence relation for the recurrence coefficients.

\section{Orthogonal polynomials on the unit circle}
In this section we will consider orthogonal polynomials on the unit circle $\mathbb{T} = \{ z \in \mathbb{C}:
|z| = 1 \}$. A very good source for the general theory is the recent set of books by Barry Simon \cite{simon}.
The sequence of polynomials $\{\varphi_n, n=0,1,2,\ldots\}$ is orthonormal on the unit circle with
respect to a weight $w$ if
\begin{equation}  \label{eq:orthounit}
   \frac{1}{2\pi} \int_0^{2\pi} \varphi_n(z) \overline{\varphi_m(z)} w(\theta)\, d\theta = \delta_{n,m}. 
\end{equation}
These polynomials are unique if we agree to make the leading coefficient positive:
\[   \varphi_n(z) = \kappa_n z^n + \cdots, \qquad \kappa_n > 0. \]
The monic polynomials are usually denoted by $\Phi_n(z) = \varphi_n(z)/\kappa_n$. An important property,
which replaces the three term recurrence relation for orthogonal polynomials on the real line,
is the Szeg\H{o} recurrence
\begin{equation}  \label{eq:recunit}
    z\Phi_n(z) = \Phi_{n+1}(z) + \overline{\alpha_n} \Phi_n^*(z),  
\end{equation}
where $\Phi_n^*(z) = z^n \overline{\Phi}_n(1/z)$ is the reversed polynomial and $\overline{\Phi}_n$ is the polynomial 
$\Phi_n$ but with complex conjugated coefficients. In \cite{simon} the recurrence coefficients $\alpha_n$ $(n=0,1,2,3,\ldots)$ are called Verblunsky coefficients. They are given by $\alpha_n = - \overline{\Phi_{n+1}(0)}$ and they
satisfy $|\alpha_n| < 1$ for $n \geq 0$ and $\alpha_{-1}=-1$.
An important relation between $\kappa_n$ and $\alpha_n$ was found by
Szeg\H{o}:
\[  \kappa_n^2 = \sum_{k=0}^n |\varphi_k(0)|^2 , \]
from which it follows that
\[   \kappa_{n+1}^2 - \kappa_n^2 = |\varphi_{n+1}(0)|^2 = \kappa_{n+1}^2 |\alpha_n|^2, \]
so that
\begin{equation}    \label{eq:kapalpha}  
 \frac{\kappa_n^2}{\kappa_{n+1}^2} = 1 - |\alpha_n|^2. 
\end{equation}

\subsection{Modified Bessel polynomials}
We will take a closer look at the orthogonal polynomials on the unit circle
for the weight
\[    w(\theta) = \exp(\lambda \cos\theta), \qquad \theta \in [-\pi,\pi]. \]
observe that $w(-\theta) = w(\theta)$, which implies that the Verblunsky coefficients are real.
Ismail \cite[pp.~236--239]{ismail} calls the resulting orthogonal polynomials the \textit{modified Bessel polynomials} since the
trigonometric moments of this weight are in terms of the modified Bessel function
\[   \frac{1}{2\pi} \int_{0}^{2\pi} e^{in\theta} w(\theta) \, d\theta = \frac{1}{\pi} \int_0^\pi
    \cos n\theta \exp(\lambda \cos \theta) \, d\theta = I_n(\lambda). \]
These polynomials appear in the analysis of unitary random matrices \cite{periwal,tracy,hisakado} 
and play an important role in
the asymptotic distribution of the length of the longest increasing subsequence of random permutations
\cite{baik}.     

Periwal and Shevitz \cite{periwal} found a non-linear recurrence relation for the Verblunsky
coefficients of these orthogonal polynomials (see also \cite{hisakado,tracy}). 
If $z=e^{i\theta}$ then
\[   w(\theta) =  \exp\left( \frac{\lambda}{2} (z + \frac1z) \right) := \hat{w}(z), \]
and this weight satisfies a Pearson equation of the form
\begin{equation}   \label{eq:pearsonunit}
 \hat{w}'(z) = \frac{\lambda}{2} (1-\frac1{z^2}) \hat{w}(z). 
\end{equation} 
Consider the integral
\begin{equation}  \label{eq:mainunit}
   \frac{1}{2\pi i} \int_{\mathbb{T}} \overline{\varphi_n(z)}\varphi_{n+1}(z) \hat{w}'(z)\, \frac{dz}{z},
\end{equation}
then by means of Pearson's equation (\ref{eq:pearsonunit}) we find
\begin{eqnarray*}
 \frac{1}{2\pi i} \int_{\mathbb{T}} \overline{\varphi_n(z)}\varphi_{n+1}(z) \hat{w}'(z)\, \frac{dz}{z}
  & = &\frac{\lambda}{4\pi i} \int_{\mathbb{T}} \overline{\varphi_n(z)}\varphi_{n+1}(z) (1- \frac{1}{z^2}) \hat{w}(z)
   \, \frac{dz}{z} \\
  & = &  \frac{\lambda}{4\pi} \int_{0}^{2\pi} \overline{\varphi_n(z)}\varphi_{n+1}(z) w(\theta) \, d\theta \\
  & &    -\ \frac{\lambda}{4\pi} \int_{0}^{2\pi} \overline{\varphi_n(z)z^2}\varphi_{n+1}(z) w(\theta) \, d\theta .
\end{eqnarray*}
The first integral on the right is zero because of orthogonality. For the second integral we use
the recurrence (\ref{eq:recunit}) for the orthonormal polynomials 
\begin{equation}  \label{eq:recunit2}
   z\varphi_n(z) = \frac{\kappa_n}{\kappa_{n+1}} \varphi_{n+1}(z) + \overline{\alpha_n} \varphi_n^*(z) 
\end{equation}
to find
\begin{eqnarray*}
 \frac{\lambda}{4\pi} \int_{0}^{2\pi} \overline{z^2\varphi_n(z)}\varphi_{n+1}(z) w(\theta) \, d\theta
   &=& \frac{\kappa_n}{\kappa_{n+1}} \frac{\lambda}{4\pi} 
   \int_{0}^{2\pi} \overline{z\varphi_{n+1}(z)}\varphi_{n+1}(z) w(\theta) \, d\theta \\
  & &  +\ \alpha_n \frac{\lambda}{4\pi} \int_{0}^{2\pi} \overline{z\varphi_n^*(z)}\varphi_{n+1}(z) w(\theta) \, d\theta.
\end{eqnarray*}
If we use (\ref{eq:recunit2}) for $n+1$ and orthogonality, then
\begin{eqnarray*}
 \frac{\kappa_n}{\kappa_{n+1}} \frac{\lambda}{4\pi} 
\int_{0}^{2\pi} \overline{z\varphi_{n+1}(z)}\varphi_{n+1}(z) w(\theta) \, d\theta
  & = & \alpha_{n+1} \frac{\kappa_n}{\kappa_{n+1}} \frac{\lambda}{4\pi} 
  \int_{0}^{2\pi} \overline{\varphi_{n+1}^*(z)}\varphi_{n+1}(z) w(\theta) \, d\theta \\
 & =&  - \frac{\lambda}{2} \frac{\kappa_n}{\kappa_{n+1}} \alpha_{n+1} \overline{\alpha_n}  
\end{eqnarray*}
because
\begin{eqnarray*}
  \varphi_{n+1}^*(z) & = & \frac{\overline{\varphi_{n+1}(0)}}{\kappa_{n+1}} \varphi_{n+1}(z) + 
   \textrm{ lower degree terms} \\
    & = &  - \alpha_n \varphi_{n+1}(z) +  \textrm{ lower degree terms}.
\end{eqnarray*}
In a similar way we have
\[  z \varphi_n^*(z) = - \alpha_{n-1} \frac{\kappa_n}{\kappa_{n+1}} \varphi_{n+1}(z) +
   \textrm{ lower degree terms} \]
so that
\[ \alpha_n \frac{\lambda}{4\pi} \int_{0}^{2\pi} \overline{z\varphi_n^*(z)}\varphi_{n+1}(z) w(\theta) \, d\theta
   = - \frac{\lambda}{2} \alpha_n \overline{\alpha_{n-1}} \frac{\kappa_n}{\kappa_{n+1}} . \]
Combining all these results gives
\begin{equation}  \label{eq:mainunit1}
  \frac{1}{2\pi i} \int_{\mathbb{T}} \overline{\varphi_n(z)}\varphi_{n+1}(z) \hat{w}'(z)\, \frac{dz}{z}
  =  \frac{\lambda}{2}\frac{\kappa_n}{\kappa_{n+1}} 
(\alpha_{n+1} \overline{\alpha_n} + \alpha_n\overline{\alpha_{n-1}}).  
\end{equation}
We can compute this integral also using integration by parts, to find
\[  \frac{1}{2\pi i} \int_{\mathbb{T}} \overline{\varphi_n(z)}\varphi_{n+1}(z) \hat{w}'(z)\, \frac{dz}{z}
  = -\frac{1}{2\pi i} \int_{\mathbb{T}} [\overline{z\varphi_n(z)}\varphi_{n+1}(z)]' \hat{w}(z)\, dz . \]
We have to be a little bit careful because $\overline{z\varphi_n(z)}$ is not analytic in the complex plane,
but on the unit circle $\mathbb{T}$ we have $\overline{\varphi_n(z)} = z^{-n}\varphi_n^*(z)$ so that
\begin{eqnarray*}
\frac{1}{2\pi i} \int_{\mathbb{T}} [\overline{z\varphi_n(z)}\varphi_{n+1}(z)]' \hat{w}(z)\, dz 
  & = & \frac{1}{2\pi i} \int_{\mathbb{T}} [z^{-n-1}\varphi_n^*(z)\varphi_{n+1}(z)]' \hat{w}(z)\, dz  \\
  & = & -\frac{n+1}{2\pi i} \int_{\mathbb{T}} z^{-n-1} \varphi_n^*(z)\varphi_{n+1}(z) \hat{w}(z)\, \frac{dz}{z} \\
  & & +\ \frac{1}{2\pi i} \int_{\mathbb{T}} z^{-n}[\varphi_n^*]'(z)\varphi_{n+1}(z) \hat{w}(z)\, \frac{dz}{z} \\
  & & +\ \frac{1}{2\pi i} \int_{\mathbb{T}} z^{-n} \varphi_n^*(z)\varphi_{n+1}'(z) \hat{w}(z)\, \frac{dz}{z} .
\end{eqnarray*}
If we use the recurrence relation (\ref{eq:recunit2}) then   
\begin{eqnarray*}
 \frac{1}{2\pi i} \int_{\mathbb{T}} z^{-n-1} \varphi_n^*(z)\varphi_{n+1}(z) \hat{w}(z)\, \frac{dz}{z}
 & = &\frac{1}{2\pi} \int_{0}^{2\pi} \overline{z\varphi_n(z)}\varphi_{n+1}(z) w(\theta)\, d\theta \\
 & = & \frac{\kappa_n}{\kappa_{n+1}},
\end{eqnarray*}
by orthogonality we find
\begin{eqnarray*}
  \frac{1}{2\pi i} \int_{\mathbb{T}} z^{-n}[\varphi_n^*]'(z)\varphi_{n+1}(z) \hat{w}(z)\, \frac{dz}{z}
  & = & \frac{1}{2\pi} \int_{0}^{2\pi} \overline{z\varphi_n'(z)}\varphi_{n+1}(z) w(\theta)\, d\theta \\
  & = & 0,
\end{eqnarray*}
and if we use 
\[ \varphi_{n+1}'(z) = (n+1) \frac{\kappa_{n+1}}{\kappa_n} \varphi_{n}(z) + \textrm{ lower degree terms}, \]
then
\begin{eqnarray*}
\frac{1}{2\pi i} \int_{\mathbb{T}} z^{-n} \varphi_n^*(z)\varphi_{n+1}'(z) \hat{w}(z)\, \frac{dz}{z} 
& = & \frac{1}{2\pi} \int_{0}^{2\pi} \overline{\varphi_n(z)}\varphi_{n+1}'(z) w(\theta)\, d\theta \\
& = & (n+1) \frac{\kappa_{n+1}}{\kappa_n}. 
\end{eqnarray*}
These computations give
\begin{equation} \label{eq:mainunit2}
   \frac{1}{2\pi i} \int_{\mathbb{T}} [\overline{z\varphi_n(z)}\varphi_{n+1}(z)]' \hat{w}(z)\, dz 
  = (n+1) \frac{\kappa_{n+1}}{\kappa_n} \left( 1 - \frac{\kappa_n^2}{\kappa_{n+1}^2} \right).
\end{equation}
Now we can  combine (\ref{eq:mainunit1}) and (\ref{eq:mainunit2}) to find
\[
 -\frac{\lambda}{2}\frac{\kappa_n}{\kappa_{n+1}} 
(\alpha_{n+1} \overline{\alpha_n} + \alpha_n\overline{\alpha_{n-1}})
 = (n+1) \frac{\kappa_{n+1}}{\kappa_n} \left( 1 - \frac{\kappa_n^2}{\kappa_{n+1}^2} \right)
\]
which, together with (\ref{eq:kapalpha}) gives
\begin{equation*}  
   -\frac{\lambda}{2} ( 1 - |\alpha_n|^2) 
(\alpha_{n+1} \overline{\alpha_n} + \alpha_n\overline{\alpha_{n-1}})
 = (n+1) |\alpha_n|^2 . 
\end{equation*}
Recall that $w(-\theta)=w(\theta)$ implies that the $\alpha_n$ are real. Hence when $\alpha_n \neq 0$ then
\begin{equation} \label{eq:uniteq}
 -\frac{\lambda}{2} ( 1 - \alpha_n^2) 
(\alpha_{n+1} + \alpha_{n-1})
 = (n+1) \alpha_n .
\end{equation}
This non-linear recurrence relation corresponds to the discrete Painlev\'e equation  d-P$_{\rm II}$
\[  x_{n+1}+x_{n-1}=\frac{x_n (\alpha n+\beta)+\gamma}{1-x_n^2} \]
with $\alpha_n=x_n$, $\alpha=\beta= -2/\lambda$ and $\gamma=0$. The initial values are
\[  \alpha_{-1} = -1, \qquad \alpha_0 = \frac{\int_{0}^{2\pi} z w(\theta)\, d\theta}{\int_{0}^{2\pi}  w(\theta)\, d\theta} =
  \frac{I_1(\lambda)}{I_0(\lambda)}. \]

\section{Discrete orthogonal polynomials}
In this section we will study certain discrete orthogonal polynomials on the integers $\mathbb{N}$. 
The orthonormality
now becomes
\begin{equation}  \label{eq:dorthonormal}
  \sum_{k=0}^\infty p_n(k)p_m(k) w_k = \delta_{n,m}, \qquad n,m \geq 0.
\end{equation}
Instead of the differential operator we will now be using difference operators, namely the forward
difference $\Delta$ and the backward difference $\nabla$ for which
\[  \Delta f(x) = f(x+1)-f(x), \qquad \nabla f(x) = f(x)-f(x-1). \]
We now have two sequences $\{a_n: n=1,2,\ldots\}$ and $\{b_n: n=0,1,2,\ldots\}$ of recurrence coefficients,
and we need two recurrence relations to determine all $a_n$ and $b_n$.

\subsection{Charlier polynomials}
Charlier polynomials are the orthonormal polynomials for the Poisson distribution
\[   w_k = \frac{a^k}{k!}, \qquad k \in \mathbb{N}, \ a > 0.  \]
Observe that 
\begin{equation} \label{eq:PearsonCharlier}  
 w_{k-1} = \frac{k}{a}\ w_k  
\end{equation}
which is the (discrete) Pearson equation for the Poisson distribution. It can also be written as 
$a \nabla w_k = (a-k) w_k$. The Pearson equation gives the following structure relation for
Charlier polynomials.
\begin{lemma}
For the orthonormal Charlier polynomials one has
\begin{equation}  \label{eq:structure}
     p_n(x+1) = p_n(x) + \frac{a_n}{a} p_{n-1}(x),
\end{equation}
where $a_n$ are the coefficients in the recurrence relation (\ref{eq:3tr}).
\end{lemma}
\begin{proof}
If we expand $p_n(x+1)$ into a Fourier series, then
\[   p_n(x+1) = \sum_{j=0}^n A_{n,j} p_j(x), \]
and if we compare the leading coefficients then $A_{n,n} = 1$. The other Fourier coefficients
are given by
\[   A_{n,j} = \sum_{k=0}^\infty p_n(k+1)p_j(k) w_k = \sum_{k=1}^\infty p_n(k)p_j(k-1)w_{k-1}, \]
and if we use (\ref{eq:PearsonCharlier}) then this gives
\[   A_{n,j} = \frac{1}{a} \sum_{k=0}^\infty k p_n(k) p_j(k-1) w_k . \]
The polynomial $xp_j(x-1)$ has degree $j+1$, hence by orthogonality $A_{n,j} = 0$ whenever $j < n-1$.
For $j=n-1$ we have 
\[  xp_{n-1}(x-1) = \frac{\gamma_{n-1}}{\gamma_n} p_n(x) + \textrm{ lower degree terms}  \]
so that (\ref{eq:angamma}) gives the desired result.
\end{proof}

Note that we can write (\ref{eq:structure}) also as
\[   \Delta p_n(x) = \frac{a_n}{a} p_{n-1}(x). \]
If we compare the leading coefficient in the latter, then $n \gamma_n = a_n \gamma_{n-1}/a$, so that
\[   a_n^2 = an. \]
We will now compute the sum
\[   \sum_{k=0}^\infty p_n^2(k+1) w_k \]
in two different ways. If we use the Pearson equation (\ref{eq:PearsonCharlier}) then
\[   \sum_{k=0}^\infty p_n^2(k+1) w_k = \sum_{k=1}^\infty p_n^2(k) w_{k-1} =
    \frac{1}{a} \sum_{k=0}^\infty k p_n^2(k) w_k = \frac{b_n}{a}. \]
On the other hand, the structure relation (\ref{eq:structure}) gives
\[  \sum_{k=0}^\infty p_n^2(k+1) w_k = \sum_{k=0}^\infty \left( p_n(k) + \frac{a_n}{a} p_{n-1}(k) \right)^2 w_k
    = 1 + \frac{a_n^2}{a^2}. \]
Combining both computations gives
\[   b_n = a + \frac{a_n^2}{a} = n+a. \]
These simple computations show that the recurrence coefficients for Charlier polynomials are  given by
\[    a_n = \sqrt{an}, \qquad b_n = n+a . \] 
  
\subsection{Generalized Charlier polynomials}
If we take the weights
\[   w_k = \frac{a^k}{(k!)^N}, \qquad k \in \mathbb{N},\ a >0 , \]
with $N \in \{1,2,3,\ldots\}$, then for $N=1$ one finds the Charlier polynomials and for $N \geq 2$
the generalized Charlier polynomials. These were introduced by Hounkonnou et al.\ in \cite{ronveaux}. 
The Pearson equation is 
\begin{equation}  \label{eq:Pearsongen}
    w_{k-1} = \frac{k^N}{a}\ w_k,
\end{equation}
which can also be written as $a\nabla w_k = (a-k^N) w_k$. For $N \geq 2$ the factor $a-k^N$ is a polynomial
in $k$ of degree greater than one, and hence the weight is no longer classical but semi-classical.
We will investigate the case $N=2$ in more detail.

\begin{lemma}  \label{lem:structure}
For $N=2$ the generalized Charlier polynomials satisfy the structure relation
\begin{equation}  \label{eq:structuregen}
   p_n(x+1) = p_n(x) + \frac{n}{a_n} p_{n-1}(x) + \frac{a_na_{n-1}}{a} p_{n-2}(x),
\end{equation}
where $a_n$ are the recurrence coefficients in the three-term recurrence relation (\ref{eq:3tr}).
\end{lemma}

\begin{proof}
If we expand $p_n(x+1)$ into a Fourier series, then
\[  p_n(x+1) = \sum_{j=0}^n A_{n,j} p_j(x). \]
Comparing coefficients of $x^n$ gives $A_{n,n} = 1$, and comparing coefficients of $x^{n-1}$ gives
$A_{n,n-1}=n/a_n$. The remaining Fourier coefficients are given by
\[  A_{n,j} = \sum_{k=0}^\infty p_n(k+1)p_j(k) w_k = \sum_{k=1}^\infty p_n(k)p_j(k-1) w_{k-1}. \]
If we use the Pearson equation (\ref{eq:Pearsongen}) then
\[  A_{n,j} = \frac{1}{a} \sum_{k=0}^\infty p_n(k)p_j(k-1) k^2 w_k. \]
The polynomial $x^2p_j(x-1)$ is of degree $j+2$ and hence by orthogonality $A_{n,j}=0$ for $j < n-2$. 
For $j=n-2$ we have
\[   x^2p_{n-2}(x-1) = \frac{\gamma_{n-2}}{\gamma_n} p_n(x) + \textrm{ lower degree terms} \]
so that
\[   A_{n,n-2} = \frac{1}{a} \frac{\gamma_{n-2}}{\gamma_n} = \frac{a_na_{n-1}}{a} , \]
where we used (\ref{eq:angamma}), which gives the required result. 
\end{proof}

The structure relation (\ref{eq:structuregen}) can also be written as
\[   \Delta p_n(x) = \frac{n}{a_n} p_{n-1}(x) + \frac{a_na_{n-1}}{a} p_{n-2}(x). \]
If we compare coefficients of $x^{n-2}$, where we use
\[   p_n(x) = \gamma_n x^n + \delta_n x^{n-1} + \cdots, \]
then we find
\begin{equation}  \label{eq:compare}
    \binom{n}{2}\gamma_n + (n-1) \delta_n  = \frac{n}{a_n} \delta_{n-1} + \frac{a_na_{n-1}}{a} \gamma_{n-2} . 
\end{equation}
If $x_{1,n} < x_{2,n} < \cdots < x_{n,n}$ are the zeros of $p_n$, then by Vi\`ete's symmetric formulas we have
\[    \frac{\delta_n}{\gamma_n} = - \sum_{k=1}^n x_{k,n} . \]
The zeros of $p_n$ are equal to the eigenvalues of the truncated Jacobi matrix
\[  J_n = \begin{pmatrix}
           b_0\ & a_1\ & & & &\\
           a_1\ & b_1\ & a_2\ & & &\\
               & a_2\ & b_2\ & a_3\ & &\\
               &     & a_3\ & b_3\ & \ddots &\\
	       &     &     & \ddots & \ddots & a_{n-1} \\
               &     &      &       & a_{n-1} & b_{n-1}
          \end{pmatrix}  \]
and the sum of all eigenvalues is the trace of the matrix, hence
\[  \frac{\delta_n}{\gamma_n} = - \sum_{j=0}^{n-1} b_j .  \]
If we use this in (\ref{eq:compare}), then
\begin{equation}  \label{eq:ba}
   \binom{n}{2} - n b_{n-1} + \sum_{j=0}^{n-1} b_j = \frac{a_n^2a_{n-1}^2}{a}. 
\end{equation}
In order to get rid of the non-homogeneous terms, we put $b_n = n+\hat{b}_n$, and the
relation becomes
\[   -n \hat{b}_{n-1} + \sum_{j=0}^{n-1} \hat{b}_j = \frac{a_n^2a_{n-1}^2}{a}. \]
Differencing both sides gives
\begin{equation}  \label{eq:Freud1}
    -na(\hat{b}_n-\hat{b}_{n-1}) = a_n^2(a_{n+1}^2-a_{n-1}^2), 
\end{equation}
which may be considered as the first Freud equation for the recurrence coefficients.

Next, we will compute
\[  \sum_{k=0}^\infty p_n(k+1)p_{n-1}(k+1) w_k  \]
in two different ways. First we use the Pearson equation (\ref{eq:Pearsongen}) to find
\begin{eqnarray*}
  \sum_{k=0}^\infty p_n(k+1)p_{n-1}(k+1) w_k & = & \sum_{k=1}^\infty p_n(k)p_{n-1}(k) w_{k-1} \\
                       & = & \frac{1}{a} \sum_{k=0}^\infty k^2 p_n(k)p_{n-1}(k) w_k  \\
                       & = & \frac{1}{a} (J^2)_{n,n-1} .
\end{eqnarray*}
The entry $(J^2)_{n,n-1}$ can be computed easily by repeatedly using the recurrence relation (\ref{eq:3tr}) and
is equal to $(J)^2_{n,n-1} = a_n(b_n+b_{n-1})$, so that
\begin{equation}  \label{eq:msum1}
   \sum_{k=0}^\infty p_n(k+1)p_{n-1}(k+1) w_k = \frac{a_n(b_n+b_{n-1})}{a} .
\end{equation}
On the other hand, we can use the structure relation (\ref{eq:structuregen}) to find
\begin{eqnarray}
 \sum_{k=0}^\infty p_n(k+1)p_{n-1}(k+1) w_k & = & \sum_{k=0}^\infty \left( p_n(k) + \frac{n}{a_n} p_{n-1}(k)
   + \frac{a_na_{n-1}}{a} p_{n-2}(k) \right) \nonumber \\
  & & \times \left(p_{n-1}(k) + \frac{n-1}{a_{n-1}} p_{n-2}(k) + \frac{a_{n-1}a_{n-2}}{a} p_{n-3}(k) \right) w_k
  \nonumber \\
  & = & \frac{n}{a_n} + \frac{(n-1)a_n}{a},  \label{eq:msum2}
\end{eqnarray} 
where the last equality follows from the orthonormality (\ref{eq:dorthonormal}). Combining
(\ref{eq:msum1}) and (\ref{eq:msum2}), and recalling that $b_n = n+\hat{b}_n$, then gives the second Freud equation
\begin{equation}  \label{eq:Freud2}
    a_n^2(\hat{b}_n+\hat{b}_{n-1}+n) = na .
\end{equation}

If we eliminate $na$ from the two equation (\ref{eq:Freud1}) and (\ref{eq:Freud2}), then
\[  -(\hat{b}_n+\hat{b}_{n-1}+n)(\hat{b}_n-\hat{b}_{n-1}) = a_{n+1}^2-a_{n-1}^2. \]
Summing both sides of this equation gives
\begin{equation}  \label{eq:almost1}
   -\hat{b}_n^2 + \sum_{k=0}^{n-1} \hat{b}_k - n \hat{b}_n = a_{n+1}^2 + a_{n}^2 - a. 
\end{equation}
Summing both sides of (\ref{eq:Freud1}) gives
\begin{equation}  \label{eq:almost2}
   a \left( \sum_{k=0}^{n-1} \hat{b}_k - n\hat{b}_n \right) = a_n^2 a_{n+1}^2 .
\end{equation}
Combining (\ref{eq:almost1}) and (\ref{eq:almost2}) then gives
\[    a_n^2a_{n+1}^2 - a(a_{n+1}^2+a_n^2) + a^2 = a \hat{b}_n^2, \]
which is equivalent to
\begin{equation}  \label{eq:bnprod}
   a\hat{b}_n^2 = (a_{n+1}^2-a)(a_n^2-a).
\end{equation}
This means that $a_n^2-a$ and $a_{n+1}^2-a$ have the same sign, and since $a_0 = 0$ we must conclude that
$a_n^2-a < 0$ for $n\geq 0$. We may therefore write
\begin{equation}  \label{eq:ancn}
    a_n^2 = a (1-c_n^2), 
\end{equation}
with $c_0=1$, and then (\ref{eq:bnprod}) becomes
\begin{equation}  \label{eq:hatbn}
    \hat{b}_n = \sqrt{a} c_n c_{n+1}.
\end{equation}
The second Freud equation (\ref{eq:Freud2}) becomes
\begin{equation}  \label{eq:cnrec}
   (1-c_n^2) \sqrt{a} (c_{n+1} + c_{n-1}) = n c_n.
\end{equation}
If we compute the coefficients $c_n$ from the recurrence relation (\ref{eq:cnrec}), then we obtain
the recurrence coefficients $b_n = n+ \hat{b}_n$ from (\ref{eq:hatbn}) and the $a_n$ from
(\ref{eq:ancn}). The non-linear recurrence relation (\ref{eq:cnrec}) corresponds to the discrete
Painlev\'e II equation 
\[    x_{n+1}+x_{n-1} = \frac{x_n(\alpha n+\beta)+\gamma}{1-x_n^2} \]
with $c_n=x_n$ and $\alpha=1/\sqrt{a}$, $\beta=\gamma=0$. We need to find the solution with $c_0=1$ and
$c_1^2 = 1-a_1^2/a$. Observe that if we require that $p_1(x)=(x-b_0)p_0/a_1$ is orthogonal to
$p_0$, then
\[     \sum_{k=0}^\infty (k-b_0) w_k  =  0  \]
so that
\[   b_0 = \frac{\displaystyle \sum_{k=0}^{\infty} k \frac{a^k}{(k!)^2}}{\displaystyle \sum_{k=0}^\infty \frac{a^k}{(k!)^2} } 
= \frac{\sqrt{a} I_1(2\sqrt{a})}{I_0(2\sqrt{a})} , \]
where 
\[     I_\nu(z) = \sum_{k=0}^\infty \frac{(z/2)^{2k+\nu}}{k! \Gamma(k+\nu+1)}  \]
is the modified Bessel function. From (\ref{eq:hatbn}) we then see that 
\begin{equation}  \label{eq:c1}
 c_1 = \frac{I_1(2\sqrt{a})}{I_0(2\sqrt{a})} . 
\end{equation} 

\begin{figure}[ht]
\centering
\rotatebox{-90}{\resizebox{3.5in}{!}{\includegraphics{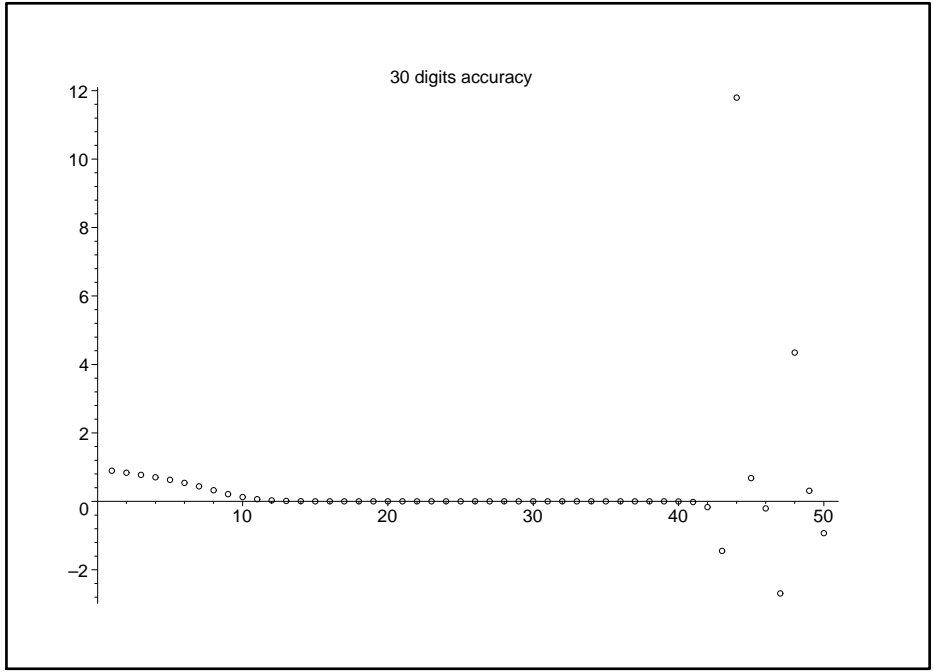}}}
\caption{The result of computing $c_n$ from d-P$_{\textrm{II}}$ using 30 significant digits}
\label{fig:dPII}
\end{figure}

The non-linear recurrence relation (\ref{eq:cnrec}) with initial conditions $c_0=1$ and
$c_1=I_1(2\sqrt{a})/I_0(2\sqrt{a})$ is again very unstable for computing all the $c_n$ recursively. One
can show that the discrete Painlev\'e equation with $\gamma=0$ and $c_0=\pm 1$ has only one solution
for which $-1 < c_n < 1$ for all $n \geq 1$ (see \cite{wva}), and this is the solution that we need since 
$a_n^2 = a (1-c_n^2)$ needs to be positive. Hence a slight deviation from the actual initial value
$c_1$ will destroy the positivity of the $a_n^2$ eventually. In Figure \ref{fig:dPII} we have
plotted the $c_n$ obtained from the recurrence relation with an accuracy of 30 digits. The $c_n$ converge
quickly to zero, but for $n$ near 40 we see that the $c_n$ deviate quite a lot from zero.

The discrete Painlev\'e II equation also arose in Section 3 for the Verblunsky coefficients of certain orthogonal
polynomials on the unit circle. Verblunsky coefficients always have the property that $|\alpha_n| < 1$ for $n \geq 0$,
hence in that case one also requires the unique solution with $\alpha_{-1}=-1$ for which $-1 < \alpha_n < 1$ for
$n \geq 0$. Observe that there is a shift in the index since we are using Verblunsky coefficients, in which case
the recurrence starts with $\alpha_{-1}=-1$.

Obviously the equation (\ref{eq:cnrec}) satisfies the discrete Painlev\'e property. Indeed, 
we have
\[   c_{n+1} = \frac{n c_n/\sqrt{a}}{1-c_n^2}-c_{n-1}, \] 
hence a singularity
will appear in $c_{n+1}$ whenever $c_n=\pm 1$. A careful analysis gives that for $c_n$ near $1$ 
\begin{eqnarray*}
    c_n & = & 1 + \epsilon \\[8pt]
    c_{n+1} & = & - \frac{n}{2\sqrt{a}}\ \frac{1}{\epsilon} - \frac{n}{4\sqrt{a}} - c_{n-1} + \O(\epsilon) \\
    c_{n+2} & = & -1 + \frac{n+2}{n} \ \epsilon + \O(\epsilon^2) \\
    c_{n+3} & = & \frac{n+1}{\sqrt{a}(n+2)} - \frac{n}{n+2} c_{n-1} + \O(\epsilon) ,
\end{eqnarray*}
and near $-1$
\begin{eqnarray*}
    c_n & = & -1 + \epsilon \\[8pt]
    c_{n+1} & = & - \frac{n}{2\sqrt{a}}\ \frac{1}{\epsilon} + \frac{n}{4\sqrt{a}} - c_{n-1} + \O(\epsilon) \\
    c_{n+2} & = & 1 + \frac{n+2}{n} \ \epsilon + \O(\epsilon^2) \\
    c_{n+3} & = & -\frac{n+1}{\sqrt{a}(n+2)} - \frac{n}{n+2} c_{n-1} + \O(\epsilon) ,
\end{eqnarray*}
so that in both cases the singularity is confined to $c_{n+1}$ and $c_{n+2}$. Observe that the critical
value $1$ for $c_n$ results in the critical value $-1$ for $c_{n+2}$, and that the critical value $-1$ for
$c_n$ results in the critical value $1$ for $c_{n+2}$. 

% added January 2006
%
\section{$q$-Orthogonal polynomials}
Here we consider orthogonal polynomials on the exponential lattice
$\{\pm q^n, n \in \mathbb{N}\}$, where $0 < q < 1$. The orthogonality is of the form
\begin{equation}  \label{eq:qortho}
  \int_{-1}^1 p_n(x)p_m(x)w(x)\, d_qx = \delta_{m,n}, 
\end{equation}
where the $q$-integral is defined by
\[  \int_{-1}^1 f(x)\, d_qx = (1-q) \sum_{k=0}^\infty f(q^k)q^k + (1-q) \sum_{k=0}^\infty f(-q^k)q^k. \]
We will only consider even weights for which $w(-x)=w(x)$, in which case the orthogonal polynomials
have the symmetry property $p_n(-x)=(-1)^n p_n(x)$, i.e., the polynomials are even when $n$ is even and odd when
$n$ is odd. The recurrence relation will then be of the form
\begin{equation} \label{eq:qrecur}
   xp_n(x) = a_{n+1}p_{n+1}(x) + a_n p_{n-1}(x), 
\end{equation}
with $p_{-1}=0$. The results in this section were obtained for the first time by Nijhoff \cite{nijhoff},
but we take a slightly different approach.

\subsection{Discrete $q$-Hermite I polynomials}
The orthonormal discrete $q$-Hermite I polynomials \cite[\S 3.28]{koekoek} are given by
\[  \int_{-1}^1 p_n(x)p_m(x) (qx;q)_\infty (-qx;q)_\infty \, d_qx = \delta_{m,n}, \]
where
\[   (x;q)_\infty = \prod_{k=0}^\infty (1-xq^k) . \]
Observe that
\[  (qx;q)_\infty (-qx;q)_\infty = \prod_{k=0}^\infty (1-x^2q^{2k}) = (x^2q^2;q^2)_\infty, \]
so that the weight can be defined as $w(x) = (x^2q^2;q^2)_\infty$. In terms of the $q$-exponential function
$E_q(z) = (-z;q)_\infty$ we have $w(x) = E_{q^2}(-x^2q^2)$, and since $E_q((1-q)z) \to \exp(z)$ when $q \to 1$
it follows that $w(\sqrt{1-q^2} x) \to \exp(-x^2)$ when $q \to 1$,
which shows that this weight is a $q$-analog of the Hermite weight. One easily finds that
\begin{equation}  \label{eq:qPearson}
     (1-x^2) w(x) = w(x/q),
\end{equation}
which is the Pearson equation for this weight on the $q$-lattice. The structure relation for the
corresponding orthogonal polynomials is in terms of the $q$-difference operator $D_q$ for which
\[   D_qf(x) = \begin{cases}\displaystyle  \frac{f(qx)-f(x)}{x(q-1)}, & \textrm{ if $x \neq 0$}, \\
                              f'(0), & \textrm{ if $x=0$}.
               \end{cases}  \]
\begin{lemma}  \label{lem:qI}
The discrete $q$-Hermite I polynomials satisfy
\begin{equation}  \label{eq:DqI}
  D_q p_n(x) = \frac{a_n}{q^{n-1}(1-q)}\ p_{n-1}.
\end{equation}
\end{lemma} 
\begin{proof}
Clearly $D_qp_n(x)$ is a polynomial of degree $n-1$ and $D_qp_n(-x)=(-1)^{n-1} D_qp_n(x)$. If we expand
this polynomial into a Fourier series, then
\[  D_qp_n(x) = \sum_{j=0}^{n-1} a_{j,n} p_j(x), \]
with 
\[  a_{j,n} = \int_{-1}^1 D_qp_n(x)\ p_j(x) w(x)\, d_qx. \]
The symmetry shows that $a_{j,n}=0$ whenever $n-j$ is even. When $n-j$ is odd then
\begin{eqnarray*}
  a_{j,n} &=& 2 (1-q) \sum_{k=0}^\infty D_qp_n(q^k) \ p_j(q^k) q^k w(q^k) \\
          &=& - 2 \sum_{k=0}^\infty [ p_n(q^{k+1}) -p_n(q^k)] p_j(q^k) w(q^k) \\
          &=& -2 \sum_{k=0}^\infty p_n(q^{k+1})p_j(q^k)w(q^k) 
            + 2 \sum_{k=0}^\infty p_n(q^k)p_j(q^k) w(q^k). 
\end{eqnarray*}
Both sums are finite since either $p_n$ or $p_j$ is an odd polynomial. Using the Pearson equation
(\ref{eq:qPearson}), and a shift in the summation index in the first sum, gives
\begin{eqnarray*}
   a_{j,n} &=& -2 \sum_{k=0}^\infty p_n(q^k)p_j(q^{k-1}) w(q^k) (1-q^{2k})
               + 2 \sum_{k=0}^\infty p_n(q^k)p_j(q^k)w(q^k) \\
            &=& \int_{-1}^1 p_n(x) \frac{p_j(x/q)-p_j(x)}{x(q-1)} w(x)\, d_qx 
                + \frac{1}{1-q} \int_{-1}^1 xp_n(x)p_j(x/q)w(x)\, d_qx.
\end{eqnarray*}
The first integral on the right is zero because of orthogonality. The second integral only gives
a contribution when $j=n-1$, in which case
\[  a_{n-1,n} = \frac{1}{1-q} \int_{-1}^1 xp_n(x)p_{n-1}(x/q) w(x)\, d_qx . \]
The recurrence relation (\ref{eq:qrecur}) gives
\[   xp_{n-1}(x/q) = qa_np_n(x/q)+qa_{n-1}p_{n-2}(x/q), \]
and since
\[     p_n(x/q) = q^{-n} p_n(x) + \textrm{ lower degree terms} \]
this gives
\[  a_{n-1,n} = \frac{q a_n}{1-q} \int_{-1}^1 p_n(x) p_n(x/q) w(x)\, d_qx = \frac{a_n}{q^{n-1}(1-q)}, \]
which gives the desired structure relation.   
\end{proof}

If we compare the leading coefficients on both sides of (\ref{eq:DqI}), then
\[    \gamma_n \frac{1-q^n}{1-q} = \gamma_{n-1} \frac{a_n}{q^{n-1}(1-q)}, \]
so that we find
\[   a_n^2 = q^{n-1}(1-q^n), \]
which are indeed the recurrence coefficients as given in \cite[\S 3.28]{koekoek}. So for these
orthogonal polynomials the recurrence coefficients can be found immediately from the
structure relation (\ref{eq:DqI}). Observe that the $a_n^2$ tend to zero exponentially fast and that
\[   \lim_{n \to \infty} a_n^2/q^{n-1} = 1, \]
and
\[   \lim_{q \to 1} \frac{a_n^2}{1-q^2} = \frac{n}{2}, \]
and the latter are the recurrence coefficients (\ref{eq:anH}) for the Hermite polynomials ($\rho=0)$.

\subsection{Discrete $q$-Freud polynomials}
A $q$-analog of the Freud polynomials with weight $\exp(-x^4)$ can be obtained by taking the weight
$w(x)=(q^4x^4;q^4)_\infty = E_{q^4}(-q^4x^4)$ on the exponential lattice. Observe that $w(\sqrt[4]{1-q^4} x) \to
\exp(-x^4)$ as $q \to 1$. This weight satisfies
\begin{equation}  \label{eq:qPearson4}
   (1-x^4)w(x) = w(x/q),
\end{equation}
and the structure relation for these semi-classical polynomials is:
\begin{lemma} \label{lem:Dqp4}
The orthonormal polynomials for which
\[   \int_{-1}^1 p_n(x)p_m(x) (x^4q^4;q^4)_\infty\, d_qx = \delta_{m,n}, \]
satisfy
\begin{equation}   \label{eq:Dqp4}
    D_qp_n(x) = \frac{B_n}{1-q}\ p_{n-1}(x) + \frac{A_n}{1-q}\ p_{n-3}(x),
\end{equation}
with
\begin{eqnarray}
     A_n & = &  \frac{a_na_{n-1}a_{n-2}}{q^{n-3}}   \label{eq:qAn} \\
     B_n & = & \frac{a_n}{q^{n-1}} \left( \sum_{j=1}^{n+1} a_j^2 - q^2 \sum_{j=1}^{n-2} a_j^2 \right). \label{eq:qBn}
\end{eqnarray} 
\end{lemma}
\begin{proof}
Expanding $D_qp_n$ into a Fourier series gives
\[    D_qp_n(x) = \sum_{j=0}^{n-1} a_{j,n} p_j(x), \]
with
\[   a_{j,n} = \int_{-1}^1 D_qp_n(x)\ p_j(x) w(x)\, d_qx. \]
Again $a_{j,n}=0$ whenever $n-j$ is even. When $n-j$ is odd then, as in the proof of Lemma~\ref{lem:qI}, we have
\[   a_{j,n} = \int_{-1}^1 p_n(x) \frac{p_j(x/q)-p_j(x)}{x(q-1)} w(x)\, d_qx
              + \frac{1}{1-q} \int_{-1}^1 x^3 p_n(x)p_j(x/q) w(x)\, d_qx, \]
where we have now used the Pearson equation (\ref{eq:qPearson4}). Again the first integral on the right
vanishes because of orthogonality. The second integral only gives a contribution when $j=n-1$ or $j=n-3$.
For $j=n-3$ we have
\[   a_{n-3,n} = \frac{1}{1-q} \int_{-1}^1 x^3 p_n(x)p_{n-3}(x)w(x)\, d_qx, \]
and since
\[    x^3 p_{n-3}(x/q) = \frac{\gamma_{n-3}}{\gamma_n} q^{-n+3} p_n(x) + \textrm{ lower degree terms} \]
we easily find
\[   a_{n-3,n} = \frac{a_n a_{n-1} a_{n-2}}{(1-q)q^{n-3}}, \]
which gives (\ref{eq:qAn}). For $j=n-1$ we have
\[   a_{n-1,n} = \frac{1}{1-q} \int_{-1}^1 x^3 p_n(x)p_{n-1}(x)w(x)\, d_qx, \]
and if we write
\begin{equation}  \label{eq:qx3pn}
    x^3p_{n-1}(x/q) = A_{n+2} p_{n+2}(x) + B_np_n(x) + \textrm{ lower degree terms}, 
\end{equation}
then the orthonormality gives
\[   a_{n-1,n} = \frac{B_n}{1-q}. \]
If we compare coefficients of $x^n$ in (\ref{eq:qx3pn}) then
\[    \delta_{n-1}q^{-n+3} = A_{n+2} \delta_{n+2} + B_n \gamma_n, \]
where $p_n(x) = \gamma_n x^n + \delta_n x^{n-2} + \cdots$, so that using (\ref{eq:qAn}) gives
\[   B_n = \frac{a_n}{q^{n-1}} \left( q^2 \frac{\delta_{n-1}}{\gamma_{n-1}} - \frac{\delta_{n+2}}{\gamma_{n+2}} \right). \]
If we compare coefficients of $x^{n-1}$ in the recurrence relation (\ref{eq:qrecur}) then
\[   \delta_{n} = a_{n+1} \delta_{n+1} + a_n \gamma_{n-1}, \]
from which one easily finds
\[     \frac{\delta_{n+1}}{\gamma_{n+1}} - \frac{\delta_n}{\gamma_n} = -a_n^2, \]
which gives
\begin{equation}  \label{eq:delgam}
   \frac{\delta_n}{\gamma_n} = - \sum_{j=1}^{n-1} a_j^2, 
\end{equation}
and using this in the formula for $B_n$ gives the desired expression (\ref{eq:qBn}). 
\end{proof}

If we compare coefficients of $x^{n-1}$ in the structure relation (\ref{eq:Dqp4}) then
\begin{equation}  \label{eq:qBn2}
    \gamma_n (1-q^n) = \gamma_{n-1} B_n . 
\end{equation}
Comparing coefficients of $x^{n-3}$ in (\ref{eq:Dqp4}) gives
\[    \delta_n (1-q^{n-2}) = \delta_{n-1} B_n + \gamma_{n-3} A_n, \]
which together with (\ref{eq:qBn2}) gives
\[   A_n = \frac{\gamma_n}{\gamma_{n-3}} \left( \frac{\delta_n}{\gamma_n} (1-q^{n-2}) 
         - \frac{\delta_{n-1}}{\gamma_{n-1}} (1-q^n) \right). \]
Together with (\ref{eq:qAn}) and (\ref{eq:delgam}) this gives 
\begin{equation}  \label{eq:aj2sum}
   a_n^2 a_{n-1}^2a_{n-2}^2 = q^{n-3} \left( q^{n-2} (1-q^2) \sum_{j=1}^{n-2} a_j^2 - (1-q^{n-2}) a_{n-1}^2
   \right). 
\end{equation}
On the other hand, if we compare (\ref{eq:qBn2}) with (\ref{eq:qBn}) then we find
\[     \frac{1-q^n}{a_n} = \frac{a_n}{q^{n-1}} \left( \sum_{j=1}^{n+1} a_j^2 - q^2 \sum_{j=1}^{n-2} a_j^2 \right), \]
which can be written as
\begin{equation}  \label{eq:almost}
    q^{n-1}(1-q^n) = a_n^2 \left( a_{n+1}^2 + a_n^2 + a_{n-1}^2 + (1-q^2) \sum_{j=1}^{n-2} a_j^2 \right). 
\end{equation}
If we take (\ref{eq:aj2sum}) with the index $n$ raised by one, then we can find
\[  (1-q^2)\sum_{j=1}^{n-2} a_j^2 = q^{-2n+3} a_{n+1}^2a_n^2a_{n-1}^2 - (1-q^2) a_{n-1}^2 - (1-q^{-n+1}) a_n^2, \]
and if we insert this in (\ref{eq:almost}) then we find the second order non-linear equation
\begin{equation}  \label{eq:qan}
   q^{n-1}(1-q^n) = a_n^2 \left( a_{n+1}^2 + q^{-n+1}a_n^2 + q^2 a_{n-1}^2 + q^{-2n+3} a_{n+1}^2a_n^2a_{n-1}^2 \right).
\end{equation}
We claim that this equation is a $q$-deformation of the discrete Painlev\'e I equation. Indeed, if we
take $x_n = a_n^2/\sqrt{1-q^4}$ then
\[    q^{n-1} \frac{1-q^n}{1-q^4} = x_n \left( x_{n+1} + q^{-n+1} x_n + q^2 x_{n-1} + 
(1-q^4) q^{-2n+3} x_{n+1}x_nx_{n-1} \right),
\]
which for $q \to 1$ converges to 
\[   \frac{n}{4} = x_n(x_{n+1}+x_n+x_{n-1}), \]
which is the discrete Painlev\'e I equation (\ref{eq:freud4}) for Freud polynomials (with $\rho=0$). 
If we put $a_n^2 = q^{n-1} y_n$ then (\ref{eq:qan}) 
can be rewritten as
\begin{equation}  \label{eq:qyn}
   q^n(y_{n+1}y_n+1)(y_{n-1}y_n+1) = 1-y_n^2.
\end{equation}
This could therefore be called a 
$q$-discrete Painlev\'e I equation  (q-P$_{\rm I}$).  

We can easily find the asymptotic behavior as $n\to \infty$. 
First observe that from (\ref{eq:qan}) we find the upper bound
\[    q^{-n+1} a_n^4 \leq q^{n-1} (1-q^n), \]
so that $a_n^4 \leq q^{2n-2}(1-q^n)$, and $a_n$ tends to zero as $n \to \infty$. Let 
$A = \limsup_{n \to \infty} a_n^2 /q^{n-1}$, then if we take a such that $a_n^2/q^{n-1}$
converges to $A$, equation (\ref{eq:qan}) gives $A^2=1$. A similar reasoning also shows that
$B=\liminf_{n \to \infty} a_n^2/q^{n-1}$ is such that $B^2=1$. Hence we may conclude that
\[     \lim_{n \to \infty} a_n^2/q^{n-1} = 1. \]

The equation (\ref{eq:qyn}) has the singularity confinement property. Indeed, a singularity occurs
for $y_{n+1}$ whenever $y_n=0$. So if we put $y_{n}=\epsilon$, then some straightforward calculus
gives
\begin{eqnarray*}
   y_n & = & \epsilon, \\
   y_{n+1} & = & q^{-n}(1-q^n) \ \frac{1}{\epsilon} - q^{-n} y_{n-1} + \O(\epsilon), \\
   y_{n+2} & = & -q^{-n-1}(1-q^n)\ \frac{1}{\epsilon} +  y_{n-1}/q + \O(\epsilon), \\
   y_{n+3} & = & - q^{-2} \frac{1-q^{n+3}}{1-q^n}\ \epsilon + \O(\epsilon^2), \\
   y_{n+4} & = & q^2 \frac{1-q^n}{1-q^{n+3}}\ y_{n-1} + \O(\epsilon).
\end{eqnarray*}
Hence the singularity is confined to $y_{n+1},y_{n+2},y_{n+3}$.

\begin{figure}[ht]
\centering
\rotatebox{-90}{\resizebox{3.5in}{!}{\includegraphics{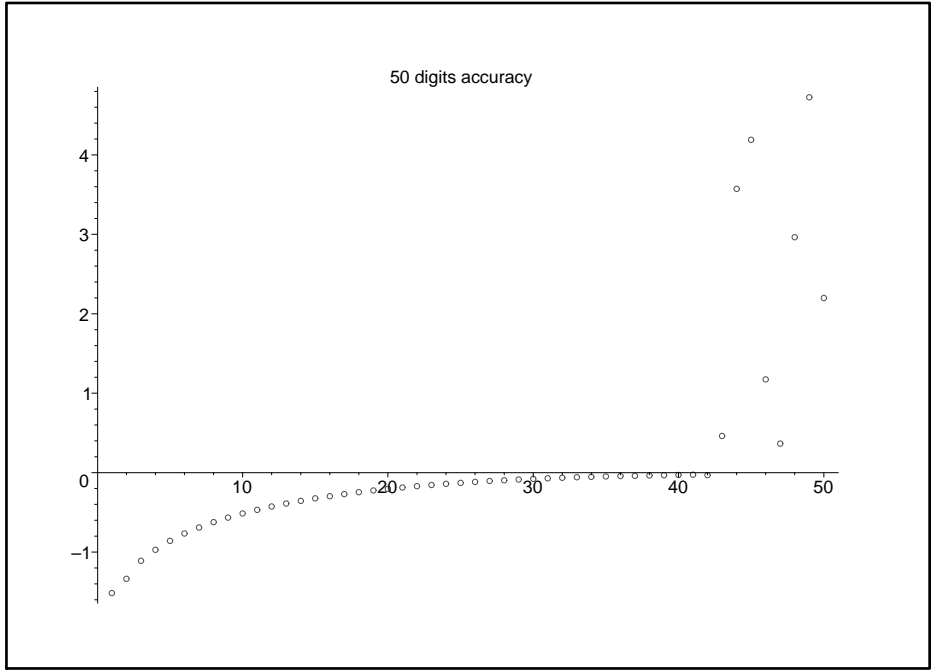}}}
\caption{The result of computing $\log |y_n|$ from (\ref{eq:qyn}) ($q=0.9$) using 50 significant digits}
\label{fig:freudy}
\end{figure}

Again the recurrence relation (\ref{eq:qan}) or (\ref{eq:qyn}) is very unstable for computing
the recurrence coefficients recursively. One can show \cite{wva} that there is again a unique solution of
(\ref{eq:qyn}) with $y_0=0$ which is positive for all $n > 0$, and this is the solution for which
$y_n=a_n^2 /q^{n-1}$. This solution is such that $y_n \to 1$ and 
\[  y_1 = a_1^2 = \frac{\int_{-1}^1 x^2 (x^4q^4;q^4)_\infty\, d_qx}{\int_{-1}^1 (x^4q^4;q^4)_\infty\,d_qx} =
  \frac{(q;q^4)_\infty}{(q^3;q^4)_\infty}, \]
where the integrals can be computed using the $q$-binomial theorem. In Figure~\ref{fig:freudy} we have computed
$\log |y_n|$ recursively for $q=0.9$ with 50 significant digits.  

\subsection{Another discrete $q$-Freud case}
If we take the weight $w(x)=(x^2q^2;q^2)_\infty (cx^2q^2;q^2)_\infty$, with $c \leq 1$, then $w$ is positive
on the $q$-lattice and it satisfies the Pearson equation
\begin{equation}  \label{eq:qPearsonc}
    (1-x^2)(1-cx^2)w(x) = w(x/q).
\end{equation}
If $c=-1+a\sqrt{1-q^4}$ then $w(\sqrt[4]{1-q^4}\, x) \to \exp(-x^4-2ax^2)$ so that this gives us a
$q$-deformation of the Freud weight $\exp(-x^4-2ax^2)$. Observe that $c=-1$ gives the discrete $q$-Freud polynomials
considered in the previous section and $c=0$ gives the discrete $q$-Hermite I polynomials.
\begin{lemma}  \label{lem:c}
The structure relation for the orthonormal $q$-polynomials with weight $w(x)=(x^2q^2;q^2)_\infty (cx^2q^2;q^2)_\infty$
on the $q$-lattice $\{\pm q^n, n \in \mathbb{N}\}$ is
\begin{equation}  \label{eq:Dqpc}
   D_qp_n(x) = \frac{\widehat{B}_n}{1-q} p_{n-1}(x) + \frac{\widehat{A}_n}{1-q} p_{n-3}(x),
\end{equation}
where
\begin{eqnarray}
     \widehat{A}_n & = & -c \frac{a_na_{n-1}a_{n-2}}{q^{n-3}}   \label{eq:qcAn} \\
     \widehat{B}_n & = & -c \frac{a_n}{q^{n-1}} \left( \sum_{j=1}^{n+1} a_j^2 - q^2 \sum_{j=1}^{n-2} a_j^2 
     - \frac{1+c}{c} \right). \label{eq:qcBn}
\end{eqnarray}  
\end{lemma}  
\begin{proof}
The proof is a straightforward copy of the proof of Lemma~\ref{lem:Dqp4}, except that one uses the
Pearson equation (\ref{eq:qPearsonc}). The Pearson equation contains the quartic polynomial
$(1-x^2)(1-cx^2)=1-(1+c)x^2+x^4$ so that one ends up with integrals of the form
\begin{eqnarray*}
  \widehat{A}_n &=& \int_{-1}^1 [(1+c)x-cx^3] p_n(x)p_{n-3}(x)w(x)\, d_qx = -c A_n\\
  \widehat{B}_n &=& \int_{-1}^1 [(1+c)x-cx^3] p_n(x)p_{n-1}(x)w(x)\, d_qx = -c B_n + (1+c) \frac{a_n}{q^{n-1}},
\end{eqnarray*} 
where $A_n$ and $B_n$ are given by (\ref{eq:qAn})--(\ref{eq:qBn}).
\end{proof}

Reasoning in the same way as in the previous section, i.e., comparing the coefficients of $x^{n-1}$ and $x^{n-3}$
in (\ref{eq:Dqp4}), one arrives at 
\begin{equation}   \label{eq:qanc}
   q^{n-1}(1-q^n) = -ca_n^2 \left( a_{n+1}^2 + q^{-n+1} a_n^2 + q^2 a_{n-1}^2 - \frac{1+c}{c} 
    - c q^{-2n+3}a_{n+1}^2a_n^2a_{n-1}^2 \right).
\end{equation} 
 If we put $a_n^2 = q^{n-1} y_n$, then this can be rewritten as
\begin{equation}  \label{eq:qync}
   (1-y_n)(1-cy_n) = q^n(cy_{n+1}y_n-1)(cy_{n-1}y_n-1),
\end{equation}
which is a more general form of the $q$-discrete Painlev\'e I equation in (\ref{eq:qyn}).

\section*{Appendix}
Several discrete Painlev\'e equations have appeared in the literature, and the list is certainly longer than
the six Painlev\'e differential equations. Sakai \cite{sakai} made a classification in terms of rational surfaces 
associated with affine root systems and the most general (elliptic) discrete Painlev\'e equation is related with the affine Weyl group symmetry of type $E_8$. Sakai's classification
does not give explicit expressions for the discrete Painlev\'e equations.
A few important discrete Painlev\'e equations are listed below. The list was compiled by Peter Clarkson
and I thank him for his permission to present it in this paper.
\def\p{Pain\-lev\'e}

\def\a{\alpha}
\def\b{\beta}
\def\c{\gamma}
\def\ga{\gamma}
\def\de{\delta}
\def\la{\lambda}
\def\q{\lambda}
\def\k{\kappa}
\def\th{\theta}
\def\ep{\varepsilon}

\def\nn{\kappa(-1)^{n+m}}
\def\anp{a_{n+1}} %\a_{n+1}
\def\an{a_{n}} %\a_{n}
\def\anm{a_{n-1}} %\a_{n-1}
\def\anpm{a_{n\pm1}} %\a_{n\pm1}

\def\bnp{b_{n+1}} %\b_{n+1}
\def\bn{b_{n}} %\b_{n
\def\bnm{b_{n-1}} %\b_{n-1}
\def\bnpm{b_{n\pm1}} %\b_{n\pm1}

\def\cnp{c_{n+1}} %\c_{n+1}
\def\cn{c_{n}} %\c_{n}
\def\cnm{c_{n-1}} %\c_{n-1}
\def\cnpm{c_{n\pm1}} %\c_{n\pm1}

\def\tfr#1#2{{\textstyle\frac{#1}{#2}}}
\def\xpp{x_{n+2}}
\def\xp{x_{n+1}}
\def\xm{x_{n-1}}
\def\xn{x_n}
\def\ypp{y_{n+2}}
\def\yp{y_{n+1}}
\def\ym{y_{n-1}}
\def\yn{y_n}
\def\z{\xi}
\def\zpp{z_{n+2}}
\def\zp{z_{n+1}}
\def\zm{z_{n-1}}
\def\zn{z_n}
\def\zhp{z_{n+1/2}}
\def\zhm{z_{n-1/2}}
\def\lpp{\la_{n+2}}
\def\lp{\la_{n+1}}
\def\lm{\la_{n-1}}
\def\ln{\la_n}

\def\dP{d-P}
\def\PI{\mbox{\rm P$_{\!\rm I}$}}
\def\PII{\mbox{\rm P$_{\!\rm II}$}}
\def\PIII{\mbox{\rm P$_{\!\rm III}$}}
\def\PIV{\mbox{\rm P$_{\!\rm IV}$}}
\def\PV{\mbox{\rm P$_{\!\rm V}$}}
\def\PVI{\mbox{\rm P$_{\!\rm VI}$}}
\def\dPI{\mbox{\rm \dP$_{\!\rm I}$}}
\def\dPII{\mbox{\rm \dP$_{\!\rm II}$}}
\def\dPIII{\mbox{\rm \dP$_{\!\rm III}$}}
\def\dPIV{\mbox{\rm \dP$_{\!\rm IV}$}}
\def\dPV{\mbox{\rm \dP$_{\!\rm V}$}}
\def\dPVI{\mbox{\rm \dP$_{\!\rm VI}$}}
\def\qPI{\mbox{\rm $q$-P$_{\!\rm I}$}}
\def\qPII{\mbox{\rm $q$-P$_{\!\rm II}$}}
\def\qPIII{\mbox{\rm $q$-P$_{\!\rm III}$}}
\def\qPIV{\mbox{\rm $q$-P$_{\!\rm IV}$}}
\def\qPV{\mbox{\rm $q$-P$_{\!\rm V}$}}
\def\qPVI{\mbox{\rm $q$-P$_{\!\rm VI}$}}
\def\aPI{\mbox{\rm a-\dP$_{\!\rm I}$}}
\def\aPII{\mbox{\rm a-\dP$_{\!\rm II}$}}
\def\aPIII{\mbox{\rm a-\dP$_{\!\rm III}$}}
\def\aPIV{\mbox{\rm a-\dP$_{\!\rm IV}$}}
\def\aPV{\mbox{\rm a-\dP$_{\!\rm V}$}}
\def\aPVI{\mbox{\rm a-\dP$_{\!\rm VI}$}}

\def\ds{\displaystyle}

\subsection*{A.1\ \ Discrete \p\ equations}

\begin{tabular}{ll}
{\dPI} & $\ds x_{n+1}+x_n+x_{n-1}=\frac{\zn +\ga (-1)^n}{x_n} + \de$\\[15pt]
{\dPII} & $\ds x_{n+1}+x_{n-1}={x_n\zn+\ga\over 1-x_n^2}$\\[15pt]
%{\dPIII} & $\ds x_{n+1}x_{n-1}=\frac{\c\de(x_n-\a q^{2n})(x_n-\b
%q^{2n})}{(x_n-\c)(x_n-\de)}$\\[15pt]
{\dPIV} & $\ds (x_{n+1}+x_{n})(x_{n}+x_{n-1}) =
\frac{(x_{n}^2-\k^2)(x_{n}^2-\mu^2)}{(x_{n}+ \zn)^2 -\ga^2}$\\[15pt]
%{\dPIV} & $\ds (x_{n+1}+x_{n})(x_{n}+x_{n-1}) =
%\frac{(x_{n}+\k+\mu)(x_{n}+\k-\mu)(x_{n}-\k+\mu)(x_{n}-\k-\mu)}
%{(x_{n}+ \zn +\ga)(x_{n} +\zn -\ga)}$\\[15pt]
{\dPV} & $\ds\frac{(x_{n+1}+x_{n}-\zp-z_{n})(x_{n}+x_{n-1}-z_{n}-\zm)}
{(x_{n+1}+x_{n})(x_{n}+x_{n-1})}=\frac{[(x_{n}-\zn)^2-\a^2][(x_{n}-\zn)^2-\b^2]}{(x_{n}-\ga^2)(x_{n}-\de^2)}$
%{\dPV} & $\ds \frac{(x_{n+1}+x_{n}-\zn-\zp)(x_{n}+x_{n-1}-\zn-\zm)}{(x_{n+1}+x_{n})(x_{n}+x_{n-1})} =
%\frac{[(x_{n}-\zn)^2-\k^2][(x_{n}-\zn)^2-\mu^2]}{(x_{n}-\ga)^2(x_{n}-\de)^2}$\\[15pt]
%{\dPV} & $\ds\frac{2\zp}{1-\xp\xn} + \frac{2z_{n}}{1-\xn\xm}
%= \mu+\nn  +\zp+\zn+ \frac{ \left[\mu-\nn\right](\ga^2-1)\xn}{(\ga+\xn)(1+\ga\xn)}$\\[10pt]
%& $\ds\qquad+ \frac{\ga(1-\xn^2)\left[\zn+\zp +
%(-1)^n(\zn-\zp-2\de) \right] }{2(\ga+\xn)(1+\ga\xn)}$
%{\dPV} & $\ds (x_{n+1}x_{n}-1)(x_{n}x_{n-1}-1)=\frac{\ga\de(x_{n}-\a)(x_{n}-1/\a)
%(x_{n}-\b)(x_{n}-1/\b)}{(x_{n}-\ga)(x_{n}-\de)}$\\[15pt]
%{\dPVI} & $\ds \frac{(x_nx_{n+1}-\zn \zp)(x_nx_{n-1}-\zn \zm)}
%{(x_nx_{n+1}-1)(x_nx_{n-1}-1)}$ \\[10pt] &\qquad\qquad
%$\ds=\frac{(x_n-\a \zn)(x_n-\zn/\a)(x_n-\b \zn)(x_n-\zn/\b)}
%{(x_n-\c)(x_n-1/\c)(x_n-\de)(x_n-1/\de)}$
\end{tabular}

\bigskip\noindent
where $\zn=\a n+\b$ and
$\a$, $\b$, $\ga$, $\de$, $\k$, $\mu$ are constants.

\subsection*{A.2\ \ $q$-discrete \p\ equations}

\begin{tabular}{ll}
{\qPII} & $\ds (x_{n+1}x_n-1)(x_nx_{n-1}-1)=\frac{\ln\lm x_n}{\a^2(x_n-\a
\ln)}$\\[15pt] {\qPII} & $\ds x_{n+1}x_{n-1}=\a 
\ln\,\frac{\ln+x_n}{x_n(x_n-1)}$\\[15pt]
{\qPIII} & $\ds x_{n+1}x_{n-1}=\frac{(x_n+\a)(x_n+\b)}{(\ga \ln
x_{n}+1)(\de \ln x_{n}+1)}$\\[15pt]
{\qPIV} & $\ds (x_{n+1}x_{n}-1)(x_{n}x_{n-1}-1)=
\frac{\ga\de(x_{n}+\a)(x_{n}+1/\a)(x_{n}+\b)(x_{n}+1/\b)} {(\ga \ln 
x_{n}+1)(\de
\ln x_{n}+1)}$\\[15pt]
{\qPV} & $\ds 
(x_{n+1}x_{n}-1)(x_{n}x_{n-1}-1)=\frac{\ga\de\q_n^{2}(x_{n}-\a)(x_{n}-1/\a)
(x_{n}-\b)(x_{n}-1/\b)}{(x_{n}-\ga\la_n)(x_{n}-\de\la_n)}$\\[15pt]
{\qPVI} & $\ds \frac{(x_nx_{n+1}-\ln\lp)(x_nx_{n-1}-\ln\lm)}
{(x_nx_{n+1}-1)(x_nx_{n-1}-1)}$ \\[15pt] &\qquad\qquad
$\ds=\frac{(x_n-\a \ln)(x_n-\ln/\a)(x_n-\b \ln)(x_n-\ln/\b)}
{(x_n-\c)(x_n-1/\c)(x_n-\de)(x_n-1/\de)}$\\[10pt]
&\qquad\mbox{with}\qquad $\ds a+b+c+d=0,\quad p+q+r+s=0$\\
\end{tabular}

\bigskip\noindent
where $\la_n=\la_0 q^n$ and
$\a$, $\b$, $\ga$ and $\de$ are constants.

\subsection*{A.3\ \ Asymmetric discrete \p\ equations}
%$$\xp = \frac{f_1(\yn)+\xn f_2(\yn)}{f_3(\yn)+\xn f_4(\yn)}, \qquad
%\ym = \frac{g_1(\xn)+\yn g_2(\xn)}{g_3(\xn)+\yn g_4(\yn)} $$

\bigskip\begin{tabular}{llll}
&$\ds\xp = \frac{f_1(\yn)+\xn f_2(\yn)}{f_3(\yn)+\xn f_4(\yn)}$
&$\ds\ym = \frac{g_1(\xn)+\yn g_2(\xn)}{g_3(\xn)+\yn g_4(\yn)}$ \\[25pt]
{$\a$-\dPI} &$\ds\xp + \xn + \yn = \de+\frac{\zn-\c}{\yn}\quad$ %$ 
\\[15pt]%&$\ds
&$\ds \yn + \ym + \xn = \de+ \frac{\zhp+\c}{\xn} $ \\[15pt]
{$\a$-\dPII} &$\ds\xp + \xn = \frac{2(\yn \zn+\c)}{1-\yn^2}\quad$ % $ 
\\[15pt]%&$\ds
&$\ds \yn + \ym = \frac{2(\xn \zhp-\de)}{1-\xn^2} $ \\[15pt]
{$\a$-\dPIII} &$\ds\xp \xn = 
\frac{(\yn-q^na)(\yn-q^nb)}{(\yn-c)(\yn-d)}\quad$ %$ \\[15pt]&$\ds
&$\ds\yn \ym = \frac{(\xn-q^n\a)(\xn-q^n\b)}{(\xn-\c)(\xn-\de)}\quad$
&$\ds\left(\frac{\a\b}{\c\de}=q\frac{ab}{cd}\right) $ \\[15pt]
\end{tabular}
\begin{tabular}{ll}
{$\a$-\dPIV} &$\ds  (\yn + 
\xn)(\xp+\yn)=\frac{(\yn-a)(\yn-b)(\yn-c)(\yn-d)}{(\yn+\c-\zn)(\yn-\c-\zn)}$ 
\\[15pt]
&$\ds  (\yn + \xn)(\xn 
+\ym)=\frac{(\xn+a)(\xn+b)(\xn+c)(\xn+d)}{(\xn+\de-\zhp)(\xn-\de-\zhp)}$ 
\\[10pt]
&\qquad\mbox{with}\qquad $\ds a+b+c+d=0$\\[15pt]
{$\a$-\dPV} &$\ds  \frac{(\yn + \xp-\zn-\zp)(\xp+\yn-\zn-\zm)}{(\yn +
\xn)(\xp+\yn)}$\\[15pt] &$\ds\qquad 
=\frac{(\yn-\zn-a)(\yn-\zn-b)(\yn-\zn-c)(\yn-\zn-d)}{(\yn-p)(\yn-q)(\yn-r)(\yn-s)}$
\\[15pt]  &$\ds  \frac{(\yn + \xp-\zn-\zhm)(\xp+\yn-\zhm-\zm)}{(\yn +
\xn)(\xn+\ym)}$\\[15pt]
&$\ds\qquad 
=\frac{(\xn-\zhm+a)(\xn-\zhm+b)(\xn-\zhm+c)(\xn-\zhm+d)}{(\xn+p)(\xn+q)(\xn+r)(\xn+s)}$\\[25pt] 
{$\a$-\qPV} &$\ds (\xn\yn-1)(\xm\yn-1) = 
q^{2n}\,\frac{(\yn-a)(\yn-b)(\yn-c)(\yn-d)}{(q^{n}-\k\yn)(q^{n}-\yn/\k)}$\\[15pt]
  &$\ds (\xn\yn-1)(\xn\yp-1) = 
q^{2n+1}\,\frac{(\xn-1/a)(\xn-1/b)(\xn-1/c)(\xn-1/d)}{(q^{n+1/2}-\mu\yn)(q^{n+1/2}-\yn/\mu)}$\\[10pt]
&\qquad\mbox{with}\qquad $\ds abcd=0$\\[25pt]
{$\a$-\dPVI} & $\ds \begin{cases} \ds x_nx_{n+1} =
\frac{\b_3\b_4(y_n-q^{n}\a_1)(y_n-q^{n}\a_2)}{(y_n-\a_3)(y_n-\a_4)}\\
\ds y_ny_{n-1} = 
\frac{\a_3\a_4(x_n-q^n\b_1)(x_n-q^n\b_2)}{(x_n-\b_3)(x_n-\b_4)} \end{cases} $
\qquad{with}\qquad$\ds\frac{\a_1\a_2}{\a_3\a_4}=q\frac{\b_1\b_2}{\b_3\b_4}$
\end{tabular}

\bigskip\noindent
where $\zn=\a n+\b$ and $a$, $b$, $c$, $d$, $p$, $q$, $r$, $s$, $\a$, 
$\b$, $\ga$ and $\de$ are constants.

\subsection*{A.4\ \ Alternative discrete \p\ equations}
\begin{tabular}{ll}%\hline
{\aPI} & $\ds x_{n+1}+x_n+x_{n-1}=\frac{\zn +\ga (-1)^n}{x_n} + \mu $\\[15pt]
& $\ds \frac{\zn}{x_{n+1}+x_{n}} + \frac{\zm}{x_{n}+x_{n-1}} =-x_n^2 
+ \c$\\[15pt]
& $\ds x_{n+1}+x_{n-1}=\frac{\zn}{x_n} + \frac{\c}{x_n^2} $\\[15pt]
& $\ds x_{n+1}+x_{n-1}=\frac{\zn}{x_n} + \c $\\[15pt]
& $\ds x_{n+1}x_{n-1}=\frac{\exp(\zn)}{x_n} + \frac{\c}{x_n^2} $\\[25pt]%\hline
{\aPII} & $\ds x_{n+1}+x_{n-1}=\frac{x_n\zn+\c}{1-x_n^2}$\\[15pt]
& $\ds \frac{\zn}{x_{n+1}x_n+1} + \frac{\zm}{x_nx_{n-1}+1}
=-x_n+\frac{1}{x_n}+\zn +\c$\\[25pt]%\hline
{\aPV} & $\ds 
\frac{(x_{n+1}+x_{n}-2\zp)(x_{n}+x_{n-1}-2\zn)}{(x_{n+1}+x_{n})(x_{n}+x_{n-1})} 
=
\frac{(x_{n}-\zhp)^2-\k^2}{(x_{n}-\ga)^2}$\\[15pt]
%&\qquad$\ds=\frac{(x_{n}-\zn-\a)(x_{n}-\zn+\a)(x_{n}-\zn-\b)
%(x_{n}-\zn+\b)}{(x_{n}-\ga)(x_{n}+\ga)(x_{n}-\de)(x_{n}+\de)}$\\[15pt]
%\hline
& $\ds  \frac{\zhp}{1-\xn\xp} + \frac{\zhm}{1-\xn\xm} = \mu +\zn +
\frac{\k\xn}{(1+\xn)^2}+\frac{1-\xn}{1+\xn}\,[\tfr12\zn + (-1)^n\ga]$\\[15pt]
\end{tabular}

\bigskip\noindent
where $ \zn=\a n+\b$ and
$\a$, $\b$, $\ga$ and $\de$ are constants.
%and $\ds\frac{\a_1\a_2}{\a_3\a_4}=q\frac{\b_1\b_2}{\b_3\b_4}$.

\subsection*{A.5\ \ Other discrete \p\ equations}
\begin{tabular}{lll}%\hline
{\dPI} &$\ds\xp + \xn = \frac{\yn \zn+\c}{\yn^2},\quad$  &$\ds
\yn + \ym = \frac{\xn \zhp+\de}{\xn^2} $ \\[15pt]
{\dPII} &$\ds\xp + \xn = \frac{\yn \zn+\c}{\yn^2-\mu^2},\quad$  &$\ds
\yn + \ym = \frac{\xn \zhp+\de}{\xn^2-\mu^2} $ \\[15pt]
{\dPIV} & $\ds  \xn\xm = 
\frac{(\yn+\zn-a)(\yn+\zn-b)}{\yn^2-\ga^2},\quad$  &$\ds \yn+\yp =
-\,\frac{\zhp+c}{\xn\ga+1}-\,\frac{\zhp+d}{\xn/\ga+1}$\\[10pt]
&\qquad\mbox{with}\qquad $\ds a+b+c+d=0$\\[15pt]
{\dPIV} & $\ds  \xn\xm = \frac{a(\yn+\zn-b)}{\yn^2-\ga^2},\quad$ 
&$\ds \yn+\yp =
\,\frac{c}{\xn}+\,\frac{\zhp+d}{\xn-1}$\\[15pt]
{\dPV} & $\ds  \xn + \xm = \frac{\zn+\mu}{1+\yn/t} + 
\frac{\zn-\mu}{1+t\yn},\quad$  &$\ds \yn\yp =
\frac{(\xn-\zhp)^2-\k^2}{\xn^2-\ga^2}$\\[15pt]
& $\ds  \xn + \xm = \frac{\ga}{1+\yn} + \frac{\zn+\de}{1-\yn},\quad$ 
&$\ds \yn\yp =
\ga\,\frac{\xn-\zn}{\xn^2-\mu^2}$\\[15pt]
\end{tabular}

\bigskip\noindent
where $\zn=\a n+\b$ and $a$, $b$, $c$, $d$, $p$, $q$, $r$, $s$, $\a$, 
$\b$, $\ga$, $\de$, $\k$ and $\mu$ are constants.

\obeylines
Katholieke Universiteit Leuven
Department of Mathematics
Celestijnenlaan 200 B
B-3001 Leuven
BELGIUM
\texttt{walter@wis.kuleuven.be}
\end{document}